\documentclass[12pt,leqno]{article} 
\usepackage{amsmath,amsfonts,amscd,ProtocoloAbeliano}

\title{\Large\bf  {The BIC of a singular foliation \\ defined by 
an abelian group of isometries}}

\author{
Martintxo Saralegi-Aranguren\thanks{UPRES-EA 2462 GŽomŽtrie-Algbre. 
FacultŽ Jean Perrin.
UniversitŽ d'Artois.   Rue Jean Souvraz SP 18.   62 307 Lens Cedex - 
France.   
{\sl saralegi@euler.univ-artois.fr}. Partially supported by the "PAI
Picasso" $n^{o} 02679SL$}
\\ {\small UniversitŽ d'Artois }
\and   
Robert Wolak\thanks{Instytut Matematyki. Uniwersytet Jagiellonski. Wl. 
Reymonta 4, 30 059
Krakow - Poland.  
{\sl robert.wolak@im.uj.edu.pl}. Partially supported by the KBN grant 5
P03A 040 20}
\\ {\small Uniwersytet Jagiellonski}
}

\begin{document}  
\maketitle

\begin{abstract}  We study the cohomology properties of the singular 
foliation $\F$ determined by an action $\Phi \colon G \times M\to M$ where the 
abelian Lie group $G$ preserves a riemannian metric  on the compact 
manifold $M$.
More precisely, we prove that the basic intersection cohomology 
$\lau{\IH}{*}{\per{p}}{\mf}$ is finite dimensional and verifies the 
PoincarŽ Duality. This duality includes two well-known      situations:
\begin{itemize}
\item[-]  PoincarŽ Duality for basic cohomology (the action $\Phi$ is  
almost free).
\item[-] PoincarŽ Duality for intersection cohomology (the group $G$ is 
compact and connected).
\end{itemize}
 \end{abstract}

\bigskip

This paper deals with an action $\Phi \colon G \times M \to M$ of an abelian
Lie group on a compact manifold $M$ preserving a riemannian metric on it. The 
orbits of this action define a singular foliation $\mathcal{F}$ on $M$. 
Putting together the orbits of the same dimension we get a stratification 
$\SF$ of $M$. This structure is still very regular. The foliation 
$\mathcal{F}$ is in fact a 
conical foliation and we can define the 
basic intersection cohomology $\lau{\IH}{*}{\per{p}}{\mf}$ (cf. \cite{SW2}).
This invariant becomes the basic cohomology $\hiru{H}{*}{\mf}$ when the 
action $\Phi$ is almost free, and the intersection cohomology 
$\lau{\IH}{*}{\per{p}}{M/G}$ when the Lie group $G$ is compact and 
connected (cf. \cite{SW2}).

The aim of this work is to prove that this cohomology  
$\lau{\IH}{*}{\per{p}}{\mf}$ is finite dimensional (cf. Theorem 
\ref{T}) and verifies the 
PoincarŽ Duality (cf. Theorem 
\ref{TPC}). This result generalizes \cite{EH}, for the almost free 
case. When $G$ is compact and connected we find these results in  \cite{GM}.

\medskip

The paper is organized as follows. In Section 1 we present the 
conical foliations. Section 2 is devoted to the study of 
isometric actions and the induced foliations. For this kind of foliations we can define the basic intersection 
cohomology. This is done in the Section 3. The main results of the work 
are proved in the Section 4.

\medskip

In the sequel $M$  is a connected, second countable, Haussdorff,   
without boundary and smooth 
(of class $C^\infty$) manifold  of dimension
$m$. 
All the maps are considered smooth unless something else is indicated.


\section{Conical foliations}  A singular foliation whose associated 
stratification is conical in the sense of Goresky-MacPherson (cf. \cite{GM}) 
is just a conical foliation.

\prg{\bf Singular foliations.}
A {\em regular  foliation}   on a manifold
$M$ is a partition $\mathcal{F}$ of $M$  by connected immersed 
submanifolds, called {\em leaves}, in such a way that each $x \in M$ 
possesses the following local model:
$$(\R^m,\mathcal{H})$$ where leaves are defined by 
$\{ dx_1 = \cdots = dx_p = 0\}$. We shall say that  
$(\R^m,\mathcal{H})$ is a {\em simple foliation}.
Notice that the leaves have the same dimension.

A {\em singular foliation}  on a manifold
$M$ is a partition $\mathcal{F}$ of $M$ by connected subsets, called {\em 
leaves},   in such a way that each $x \in M$ 
possesses the following local model:
$$
(\R^{m-n} \times \R^{n},\mathcal{H} \times
\mathcal{K})
$$ 
where
$(\R^{m-n},\mathcal{H})$ is  a simple foliation and
$(\R^{n},\mathcal{K})$ is  a singular foliation having the origin as 
an unique 0-dimensional leaf. When $n=0$ we just have a regular foliation.
The above condition is an empty one when $x$ is an isolated 
0-dimensional leaf.
Notice that the leaves are connected immersed 
submanifolds whose dimensions may vary .
This local model is exactly the local model of a foliation 
of Sussman \cite{Su} and Stefan \cite{St}; 
so, these foliations are singular  in our sense.

  Classifying the points of
$M$ following the dimension of the leaves one gets a {\em stratification}
$\SF$ of 
$M$ whose elements are called {\em strata}. The foliation is regular
 when this stratification has just
one stratum
$\{ M
\}$.

A smooth map $f \colon (M,{\mathcal  F}) \rightarrow  (M',{\mathcal  
F'})$  between singular foliated manifolds is {\em foliated} if it 
sends the leaves of ${\mathcal  F}$ into the leaves of ${\mathcal  
F'}$. When $f$ is an open foliated embedding then  it preserves the dimension of the 
leaves and therefore it sends the strata of $M$ into the strata of 
$M'$.  
 Examples, more 
properties and the singular version of the Frobenius theorem the 
reader can find in \cite{SW2,BA,St,Su,VA}.

\prgg {\bf Examples}.

\Zati In any open subset $U \subset M$ we have the singular foliation
$\mathcal{F}_U = $ \{connected components of 
 $L \cap U \ / \ L \in \mathcal{F}$ \}. The associated stratification is
$
\mathcal{S}_{\mathcal{F}_U} = \{ \hbox{connected components of } S
\cap U \ /
\ S\in \SF\}.
$

\zati  Consider $(N,\mathcal{K})$ a connected regular foliated manifold. In 
the product 
$N
\times M$ we have the singular foliation $\mathcal{K} \times
\mathcal{F} =
\{ L_1 \times L_2 \ / \ L_1 \in \mathcal{K},  L_2 \in \mathcal{F}\}$.  The 
associated stratification is
$
\SKF = \{  N \times  S \ / \  S \in 
\SF\}.
$

\zati  Consider $\S^{n-1}$ a sphere with positive dimension endowed with a singular 
foliation
$\mathcal{G}$ without 0-dimensional leaves. Identify the disk
$\D^{n}$ with the cone
$c\S^{n-1} = \S^{n-1} \times [0,1[  \ / \ 
\S^{n-1}\times
\{ 0\}$ by the map  $ x \mapsto [x/||x||,||x||]$ where  $[u,t]$ is a 
generic element of the cone. We shall consider on 
$\D^{n}$  the singular foliation
$$ c\mathcal{G} =  \{ F \times \{ t \}  \ / \  F \in \mathcal{G}, t
\in ]0,1[
\}
\cup \{ \vartheta \},
$$  where $\vartheta$ is the vertex $[u,0]$ of the cone. The induced 
stratification is
$$
\ScG =  
\{ S \times ]0,1[  \ / \  S \in \SG \}Ê
\cup \{\vartheta\}, $$ 
since $\mathcal{G}$ does not possesses $0$-dimensional leaves. 
For technical reasons we allow $n$ to take the value 0; in this case 
$\S^{n-1}=
\emptyset$ and $c\S^{n-1} = \{ \vartheta \}$.

%

Same considerations apply to the $\infty$-cone 
$\cin \S^{n-1} = \S^{n-1} \times [0,\infty[\Big/ \S^{n-1}  \times \{ 0 \}$. 
In particular, $\R^{n} = 
\cin\S^{n-1}$. Notice that the map $f \colon (c\S^{n-1},c\mathcal{G}) \to 
(\cin \S^{n-1},\cin\mathcal{G}) $, defined by $f[\theta,t] = [\theta, 
\tan(t\pi/2)]$, is a foliated diffeomorphism.

\medskip

 The strata of a singular  foliation are not necessarily manifolds  and 
their relative position can be very wild. Consider
$(\R,
\mathcal{F})$ where
$\mathcal{F}$ is given by a vector field $f
\frac{\partial}{\partial t}$; there are two kind of strata. The connected 
components of
$f^{-1}(\R-\{ 0\})$ and these of
$f^{-1}(0)$. In other words, any connected closed subset of $\R$ can be a 
stratum.
In order  to support an intersection cohomology, the 
stratification
$\mathcal{S}_{\mathcal{F}}$ asks for a certain amount of regularity and  
conicalicity (see
\cite{GM} for the case of  stratified pseudomanifolds). This leads us to 
introduce the following definition.

\prg {\bf Conical foliations}. This definition is made by induction 
on the dimension of $M$. 
A singular foliation $(M,\mathcal{F})$ is 
said to be a {\em conical foliation} if
 for any point $x \in M$  we can find a foliated diffeomorphism 
$$
\phii \colon (\R^{m-n} \times c\S^{n-1},\mathcal{H} \times c\mathcal{G}) \TO 
(U,\mathcal{F}_U),
$$
where
\begin{itemize}
\item $U \subset M$ is an open neighborhood of $x$,
\item $(\S^{n-1},\mathcal{G})$ is a conical foliation without  0-dimensional 
leaves, called {\em link} of $x$, and 
\item $\phii(0,\vartheta) = x$.
\end{itemize}
We shall say that $(U,\phii)$ 
is a {\em conical chart} of $x$.  Notice that, if $S$ is the
stratum containing $x$ then $\phii^{-1}(S \cap  U) =
\R^{m-n} \times \{ \vartheta \}$. 
A regular foliation can be considered as a 
conical foliation for which the link is empty.
If $m =1$, a conical foliation has to be regular. If $m=2$ a conical foliation
is either regular or it has leaves of dimension 0 and 1. The link of a 
singular leaf is $\sbat$ with the one leaf foliation. The typical 
example is given by the standard action of $\sbat$ on $\R^{1}$.

We also say that $(M \times [0,1[^p,\mathcal{F} \times \mathcal{I})$, where 
$\mathcal{I}$ is the pointwise foliation of  $[0,1[^p$, is a 
conical foliated manifold.

 Notice that each stratum is an embedded submanifold of $M$. 
 Although a point
$x$ may have several charts the integer $n$ is an invariant: it is the 
codimension of the stratum containing $x$. This integer cannot to be 1 
since the conical foliation $(\S^{n-1},\mathcal{G})$  has not 0-dimensional 
leaves.

The above local description implies some important facts about the 
stratification ${\cal 
S}_{\mathcal  F}$.  Notice for example that the
family of strata is finite in the compact case and locally finite 
in the general case.
The closure $\overline{S}$ of a stratum $S\in \SF$ is a union of strata. 
Put $ S_1 \preceq S_2$ if $S_1, S_2 \in \SF$ and
$S_1 \subset \overline{S_2}$. This relation is an order relation and 
therefore
$(\SF,\preceq)$ is a poset. 

The  {\em depth}  of $\SF$, written $ \depth  \SF$, is defined to be the largest $i$ 
for which there 
exists a chain of strata $S_0  \prec S_1 \prec \cdots \prec S_i$. So, 
$ \depth  \SF = 0$ if and only if the foliation 
${\mathcal  F}$ is regular. It is always finite because of the locally 
finiteness of $\SF$.  We also have $ \depth \SFU \leq  \depth   \SF$ for any 
open subset $U \subset M$ and  $ \depth  {\SG} =  \depth  {\SHG} < \depth  {\SHcG}$ (cf. 1.1.1).

The minimal strata are exactly the closed strata. An inductive 
argument shows that the maximal strata are the strata of dimension 
$m$. They are called {\em  regular strata} and the others {\em 
singular strata}.  Since the codimension of singular strata is at 
least 2, then only one regular stratum $R_\mathcal{F}$ appears, which is an open 
dense subset.  The union of singular 
strata shall be denoted by $\Sigma_\mathcal{F}$.

The  restriction of $\mathcal{F}$ to a stratum $S$ defines a regular foliation 
$\mathcal{F}_S$. Moreover, if $ S_1 \prec S_2$ then $\dim 
\mathcal{F}_{S_1} <\dim 
\mathcal{F}_{S_2}$. So, the  the biggest leaves of $\mathcal{F}$ live 
in 
$R_\mathcal{F}$. 
The {\em dimension} $\dim \mathcal{F} $ of the foliation $\mathcal{F}$ is the dimension of 
these leaves, that is $\dim \mathcal{F} = \dim \mathcal{F}_{R_\mathcal{F}}$.
Notice that for any singular stratum $S$ we have the equality
\begin{equation}
\label{1Dic}
\dim \mathcal{F} - \dim \mathcal{F}_S =  \dim \mathcal{G},
\end{equation}
where $(\S^{n-1} ,\mathcal{G})$ is the link of a point of $S$ (see 
\refp{2Dic}).
\bigskip
\prgg {\bf Examples}.

\Zati  Any regular foliation is a conical foliation.

\zati In the examples of section 1.1.1 we obtain conical foliations if we 
replace singular foliation by conical foliation.

\zati {\bf Foliated bundles}. Consider $(N,\mathcal{N})$ a conical 
foliation, $B$ a manifold, $\widetilde{B}$ its universal covering and $h 
\colon \pi_{1}(B) \to \Diff(N,\mathcal{N})$ a representation into the 
group of diffeomorphisms of $N$ preserving $\mathcal{N}$. The 
foliated bundle is the crossed product $M =\widetilde{B} \times_{\pi_1(B)} 
N$ 
endowed with the foliation $\mathcal{F}$ whose leaves are obtained from $\widetilde{B} 
\times L$, $L\in \mathcal{N}$. This foliated bundle is conical. The dynamic 
of the leaves, strata, \ldots  can be very complicated.
\begin{itemize}
    \item[1.] Put $B =\sbat$, $N = c_{\infty}\S^n$, $n>0$,  $\mathcal{N} = c 
    \mathcal{G}$, $\mathcal{G} $ the one-leaf foliation and 
$h(1)$ defined by $h(1)(u)=u/2$. The minimal stratum $S =\sbat$ 
does not possesses a small saturated tubular neighborhood.
     
\item[2.] Put $B = \sbat$, $N = c\S^{2n+1}$,  $\mathcal{G} $ the Hopf foliation 
and $h(1)$ defined by  $$h(1)[(z_{1},\ldots,z_{n}),t]= 
[(e^{2\pi i t} \cdot z_{1},\ldots,e^{2\pi i t} \cdot 
z_{n}),t]$$ and $\mathcal{F} = c 
    \mathcal{G}$. Here the 
leaves are diffeomorphic to cylinders (in irrational 
$t$-levels) or to  tori (in non-zero rational 
$t$-levels) or to a circle (for $t=0$) and they are completely mixed up. The minimal stratum $S =\sbat$ 
possesses  small saturated tubular neighborhoods, but these 
neighborhoods do not retract to $S$ through a foliated retraction.
\end{itemize}

\zati {\bf Singular Riemannian foliations}. A {\em singular Riemannian 
foliation}, SRF for short, is a singular foliated
manifold
$(M,\mathcal{F})$ endowed with a Riemannian metric  such that every 
geodesic that is perpendicular at one point to a leaf
remains perpendicular to every leaf it meets
(cf. \cite{Mo}). The restriction of $\mathcal{F}$ to each stratum $S$ of 
$\SF$ is a regular Riemannian foliation. The local 
structure of
$\mathcal{F}$ can be described as follows (cf.
\cite{Mo} and \cite{BM}). Given $x \in M$ we can find a chart modeled on
$(\R^{m-n}
\times
\R^{n},\mathcal{H} \times \mathcal{K})$ where 
$(\R^{m-n},\mathcal{H} )$ is a simple foliation,
$(\R^{n},\mathcal{K})$ is a SRF having the origin as unique 
0-dimensional leaf and being invariant by homotethies.

The sphere $\S^{n-1}$ inherits a 
singular Riemannian foliation $\mathcal{G}$ with 
$\mathcal{K}= c\mathcal{G}$ since 
$\mathcal{K}$ is invariant by homotheties \cite{Mo}. We conclude that a 
SRF is a conical foliation.

\zati {\bf Isometric actions}. 
The foliation determined by an  isometric action induces a singular 
Riemannian foliation and therefore a conical foliation. This example will 
be treated more extensively in the second Section.

\zati {\bf The product}. Given two conical foliated manifolds 
$(M_1,\mathcal{F}_1)$ and $(M_2,\mathcal{F}_2)$ their product $(M_1 \times 
M_2,\mathcal{F}_1 \times \mathcal{F}_2)$  is also conical. 

\medskip
\pro Let us see that. We proceed by induction on $\depth {{\SF}_1} + \depth 
{{\SF}_2}$. When this 
number is 0 the product foliation is conical since is regular. 
 The proof of the general case is just a local matter. Taking into account 
(b) it suffices to consider the product $(c\S^{n_1},c\mathcal{G}_1) \times 
(c\S^{n_2},c\mathcal{G}_2)$, where $(\S^{n_1},\mathcal{G}_1)$ and 
$(\S^{n_2},\mathcal{G}_2)$ are conical. By induction hypothesis  it 
suffices to check the local conicalicity near $(\vartheta_1,\vartheta_2)$.
We are going to construct a conical foliated manifold 
$(\S^{n_1+n_2+1},\mathcal{G})$ and a foliated embedding
$\zeta \colon (c\S^{n_1+n_2+1},c\mathcal{G}) \to 
(c\S^{n_1},c\mathcal{G}_1) \times (c\S^{n_2},c\mathcal{G}_2)$ sending 
$\vartheta$ to $(\vartheta_1,\vartheta_2)$.

\hspace{.5cm}  Since the foliation on the cone $c\S^{n_i}$ is tangent to 
each sphere  $S^{n_i} \times \{t\}$ then the submanifold 
$N = \{ ([\theta_1,t_1],[\theta_2,t_2]) \in c\S^{n_1} \times c\S^{n_2} \ / 
\ t_1^2 + t_2^2 = 1/2\}$ is a foliated submanifold. Moreover, the map 
$$
\psi \colon (N \times ]0,1[,(c\mathcal{G}_1 \times c\mathcal{G}_2)_N\times 
\mathcal{I}) \longrightarrow (c\S^{n_1} \times c\S^{n_2} ,c\mathcal{G}_1 
\times c\mathcal{G}_2 ) ,
$$
defined by $\psi ([\theta_1,t_1],[\theta_2,t_2],t) = ([\theta_1,t \cdot 
t_1],[\theta_2,t \cdot t_2])$ is a foliated embedding.
 We claim that $(N,(c\mathcal{G}_1 \times c\mathcal{G}_2)_N)$ is conical. 
Again it is a local question, so it suffices to prove that $(N - \{ t_1 = 
0 \}, (\mathcal{G}_1 \times \mathcal{G}_2)_{N - \{ t_1 = 0 \}})$ is 
conical (similarly for $N - \{ t_2 = 0 \})$. Since it is foliated 
diffeomorphic to $(\S^{n_1} \times ]0,1[ \times c\S^n_2, \mathcal{G}_1 
\times \mathcal{I} \times c\mathcal{G}_2)$, where $\mathcal{I}$ is the 
0-dimensional foliation of $]0,1[$,  then it suffices to apply the 
induction hypothesis and we get the claim. 

\hspace{.5cm} Notice that $N$ is diffeomorphic to the sphere 
$\S^{n_1+n_2+1}$, by $([\theta_1,t_1],[\theta_2,t_2]) \mapsto (\sqrt{2} 
t_1 \cdot \theta_1,\sqrt{2} t_2 \cdot \theta_2)$.
The induced foliation $(\S^{n_1+n_2+1},\mathcal{G})$ is conical. Under 
this diffeomorphism, the foliated embedding $\psi$ becomes the foliated 
embedding
$$
\zeta \colon (\S^{n_1+n_2+1} \times ]0,1[,\mathcal{G} \times \mathcal{I}) 
\longrightarrow (c\S^{n_1} \times c\S^{n_2} ,c\mathcal{G}_1 \times 
c\mathcal{G}_2 ) ,
$$
defined by $\zeta((u,v) = (x_0,\ldots,x_{n_1},y_0,\ldots y_{n_2}),t) = 
([\frac{u}{||u||}, \frac{t\cdot ||u||}{\sqrt{2}}],[\frac{v}{||v||}, 
\frac{t\cdot ||v||}{\sqrt{2}}])$. But $\zeta$ extends to the foliated 
embedding 
$\zeta \colon (c\S^{n_1+n_2+1} ,c\mathcal{G}) \longrightarrow (c\S^{n_1} 
\times c\S^{n_2} ,c\mathcal{G}_1 \times c\mathcal{G}_2 )
$
by putting $\zeta(\vartheta) = (\vartheta_1,\vartheta_2)$. 
This ends the proof
\qed

\medskip

The following result will be useful in the sequel.

\bL
\label{embed}
Let $(M,\mathcal{F})$ and $(M',\mathcal{F}')$ be two conical foliated 
manifolds with same dimension and let $f \colon M \times [0,1[ \to M' \times 
[0,1[$ be  an embedding. If the restriction 
\begin{equation}
\label{hyp}
 f\colon (M \times ]0,1[ , \mathcal{F} \times \mathcal{I}) 
\to (M' \times ]0,1[ , \mathcal{F}' \times \mathcal{I})
\end{equation}
is foliated then 
\begin{equation}
\label{concl}
f \colon (M \times [0,1[ , \mathcal{F} \times \mathcal{I}) 
\to (M' \times [0,1[ , \mathcal{F}' \times \mathcal{I})
\end{equation}
is also foliated.
\eL
\pro
Notice first that, when the two foliations are regular, the result comes directly from the 
local description of $f$. 

Consider  on the other hand $S\in \SF$ a minimal stratum. From \refp{hyp} there exists $S' 
\in \SFp$ with $f(S \times ]0,1[) \subset S' \times ]0,1[$ and therefore 
$f(S \times [0,1[) \subset \overline{S'} \times [0,1[$. We claim that 
$$f(S \times [0,1[) \subset S' \times [0,1[.$$ For that purpose, let us suppose that 
there exists $S'_0 \in \SFp$ with $S'_0 \prec S'$ and $f(S \times \{ 0 \}) \cap (S'_0 \times \{ 
0 \}) \not= \emptyset$.
Since $f(M \times [0,1[) $ is an open subset of $M' \times [0,1[$ 
then $f(M \times ]0,1[) \cap (S'_0 \times ]0,1[ )\not= \emptyset$.
But this is not possible since the map 
$
 f\colon (M \times ]0,1[ , \mathcal{F} \times \mathcal{I}) 
\to (f(M \times ]0,1[ ), \mathcal{F}' \times \mathcal{I})
$
is a foliated diffeomorphism and $S \times ]0,1[$ a minimal stratum of 
${\sf S}_{\mathcal{F} \times \mathcal{I}}$.

We proceed now by induction on $\depth \SF$. When this depth is 0 then 
$\mathcal{F}$ is a regular foliation and the above considerations  give 
$f(M \times [0,1[) \subset R' \times [0,1[$, where $R'$ is a regular 
stratum of $\SFp$. We get the result since the two foliations are regular.

Consider now the general case. Denote $S_{min}$ 
the union of closed strata of $\mathcal{F}$. By induction hypothesis the restriction 
$$
f \colon ((M - S_{min}) \times [0,1[, \mathcal{F} \times \mathcal{I}) \TO (M'  
\times [0,1[ , \mathcal{F}' \times \mathcal{I}) 
$$
is a foliated map. Consider now $S \in \SF
$ a singular stratum. We have seen that there exists $S'\in \SFp$ with 
$f(S \times [0,1[) \subset S' \times [0,1[$. It remains to prove that 
$$
f \colon (S \times [0,1[, \mathcal{F} \times \mathcal{I}) \TO (S'  
\times [0,1[ , \mathcal{F}' \times \mathcal{I}) 
$$
is a foliated map. We get the result since the two foliations are regular.
\qed

 \prg {\bf Local blow up}. Consider a conical chart $(U,\phii)$ of a 
 point $x$ of a singular stratum $S$ of $M$.  The composition map 
$$
P_{\phibas} \colon 
(\R^{m-n} \times  \S^{n-1} \times [0,1[ ,\mathcal{H} \times c
\mathcal{G} \times \mathcal{I})
\TO 
(U,\mathcal{F}_U),
$$
defined by $P_{\phibas}(u,\theta,t) = \phii(u,[ \theta,t])$, is said to be a {\em local blow 
up}. It is a foliated smooth map verifying
\begin{itemize}
\item[-] The restriction $P_{\phibas}\colon \R^{m-n} \times  \S^{n-1}\times 
]0,1[ \TO U-S$ is a diffeomorphism.
\item[-] The restriction $P_{\phibas}\colon \R^{m-n} \times \S^{n-1}
\times \{0\} \TO U\cap S$ is a fiber bundle whose fiber is just $\S^{n-1} $.
\end{itemize}
In other words, the local blow up replaces each point of the minimal stratum 
by a link.

The following result shows that the local blow up is essentially unique and 
that the link of $x$ is also unique.

\bP
\label{lema}
Let $(U_1,\phii_1)$, $(U_2,\phii_2)$ be two charts of a point $x$ 
of $M$  with $U_1 
\subset U_2$. 
There exists  an unique embedding 
$$
\Phi_{1,2} \colon \R^{m-n} \times  \S^{n-1} \times [0,1[ \to 
\R^{m-n} \times  \S^{n-1}\times [0,1[
$$
 making the following diagram commutative 
$$ 
\begin{picture}(170,60)(00,-5)

\put(0,0){\makebox(0,0){$ \R^{m-n} \times  c\S^{n-1}$}}
\put(180,0){\makebox(0,0){$ \R^{m-n} \times  c\S^{n-1}.$}}
\put(180,50){\makebox(0,0){$ \R^{m-n} \times  \S^{n-1}\times [0,1[$}}
\put(0,50){\makebox(0,0){$\R^{m-n} \times  \S^{n-1}\times [0,1[$}}

\put(0,40){\vector(0,-1){26}} 
\put(180,40){\vector(0,-1){26}} 
\put(55,0){\vector(1,0){66}}
\put(75,50){\vector(1,0){46}}

\put(88,12){\makebox(0,0){$\phii_2^{-1} \rond \phii_1$}}
\put(98,58){\makebox(0,0){$\Phi_{1,2}$}}
\put(-10,27){\makebox(0,0){$P$}} 
\put(190,27){\makebox(0,0){$P$}} 
\end{picture} 
$$ 
where $P(u,\theta,t) = (u,[\theta,t])$.

Moreover, if the two charts are 
conical charts, modeled respectively on
$(\R^{m-n} \times c \S^{n-1}, \mathcal{H}_{1} \times c 
\mathcal{G}_{1})$ and $(\R^{m-n} \times c \S^{n-1}, \mathcal{H}_{2} \times c 
\mathcal{G}_{2})$, then 
$$
\Phi_{1,2} \colon (\R^{m-n} \times  \S^{n-1}  \times [0,1[,\mathcal{H}_1 \times 
\mathcal{G}_1 \times 
\mathcal{I} ) \to (\R^{m-n} \times  \S^{n-1}  \times [0,1[,\mathcal{H}_2 \times 
\mathcal{G}_2   \times 
\mathcal{I})
$$
is a foliated embedding.
\eP

\pro Notice that 
$\phii_2^{-1} 
 \rond \phii_1 \colon \R^{m-n} \times  c\S^{n-1} \to \R^{m-n} \times  
 c\S^{n-1}
 $ 
 is an embedding. 

{\em Uniqueness}. It comes from density of $\R^{m-n} \times  \S^{n-1} \times 
]0,1[$ and from the fact that 
$P \colon \R^{m-n} \times  \S^{n-1} \times 
]0,1[
\to \R^{m-n} \times (c\S^{n-1}-\{ \vartheta\})$ is a diffeomorphism. 

{\em Existence}. Denote $ (f,g)\colon \R^{m-n} \times  \R^{n} \to 
\R^{m-n} \times  \R^{n}$ the embedding $\phii_2^{-1} 
 \rond \phii_1$. The components $f$ and $g$ are smooth with 
$g(0) =0$.
So, the map $h \colon \R^{m-n} \times  
\S^{n-1}  \times [0,1[ \to \R^{n}$ defined by $h(u,\theta,t) = 
g(u,t \cdot \theta)/t$  is smooth and 
without zeroes. Finally, we define $$\Phi_{1,2}(u,\theta,t) = \left(f(u,t \cdot 
\theta),\frac{h(u,\theta,t)}{|| h(u,\theta,t)||},
t \cdot || h(u,\theta,t)||\right).$$

{\em Embedding}. Notice first the following. Since $(f,g)$ is an embedding 
with $g(0)=0$ then each $g(u,-) \colon \R^{n} \to \R^{n}$ is an 
embedding. Put $G_u$  its derivative at 0, which is an isomorphism. By 
construction we have $h(u,\theta,0) = G_u (\theta)/||G_u (\theta)||$. So, 
each restriction $\Phi_{1,2} \colon \{ u \} \times \S^{n-1}\times \{ 0 \} \to 
\{ f(u,0) \} \times \S^{n-1} \times \{ 0 \}$ is a diffeomorphism. 

Now, consider two points (resp. two vectors) on $\R^{m-n} \times  \S^{n-1}
\times [0,1[$ sent by $\Phi_{1,2}$ to the same point (resp. vector).
Since $\phii_2^{-1} 
 \rond \phii_1 $  is an embedding then they live on a fiber $\{ u \} 
\times \S^{n-1}\times \{ 0 \}$. Since $G_u$ is  an isomorphism, we get the 
claim.

{\em Foliated}.  The embedding $\Phi_{1,2}$ extends the foliated map $\phii_2^{-1} 
 \rond \phii_1$. We 
 apply  Lemma \ref{embed}.
\qed

An important consequence of this Proposition is that the links 
$(\S^{n-1},\mathcal{G})$ of two points of the same stratum $S$ are foliated 
diffeomorphic. We shall write $\mathcal{G} _{S}= \mathcal{G}$.


\section{Foliations determined by an abelian isometric action.}

We deal in this paper with an abelian isometric action 
$
\Phi \colon G \times M \to M$ defined on a compact smooth manifold. 
 As we study  the induced foliation $\F$, we may assume that the group 
 $G$ is connected; it suffices to replace $G$ by the connected 
 componenet $G_{0}$ containing the unity element.

\smallskip 

For technical reasons, we also need  to work  
with non compact manifolds. The tame actions are introduced for this 
purpose.

\prg {\bf Tame actions}. A smooth action $\Phi \colon G \times M 
\to M$ of a abelian Lie group on a smooth manifold $M$ is  {\em tame} 
if it extends to a smooth action $\Phi \colon K \times M 
\to M$ where $K$ is an abelian compact Lie group containing $G$. 
We say that $K$ 
is a {\em tamer group} of $G$.
When $M$ is compact, this notion is equivalent to that of isometric 
action (cf.  \cite{K}).

When necessary, we 
can suppose that  the 
group $G$ is dense on $K$. 

The restriction of  the action of $G$ to a $K$-invariant submanifold of $M$ is again a 
tame action. For a subgroup $H\subset G$ the restriction $\Phi \colon H 
\times M \to M$ is also a tame action.

 The connected components of the orbits of the tame action $\Phi$  determine 
 a foliation $\mathcal{F}$ 
 on $M$. Since the action $\Phi \colon G \times M \to M$ is isometric 
 then $\F$ is a conical foliation (cf. 1.2.1 (e)). 


\prgg {\bf Three particular actions}.

In this work we shall use the particular actions we describe now.

\Zati {\em The action $\Xi \colon G \times K \to K$, defined by  $\Xi(g,k) = g \cdot 
k$. }

\nt 
The action is tame since it extends to the action  $\Xi \colon K \times K \to K$ 
defined by $\Xi(g,k) = g\cdot k$. We write $\F_{K}$ the induced foliation,
we have $\dim \F_{K} = \dim G$. 
For each $u \in \mathfrak{g}$, Lie algebra of $G$,  we  shall write $X_u \in 
\mathcal{X}(K)$ the associated fundamental vector field.

\zati {\em 
The action $\Psi \colon G \times K/H \to K/H$, defined 
by  $\Psi(g,kH) = (g \cdot 
k)H$, where $H \subset K$ a closed subgroup. }

\nt 
The action is tame since it extends to the action $\Psi \colon K \times K/H 
\to K/H$ 
defined by $\Psi(g,kH) = (g \cdot 
k)H$.
We write $\F_{K/H}$ the induced foliation, we have $\dim \F_{K/H} = \dim G H/ H = \dim G/(G \cap H)$. 
For each $u \in \mathfrak{g}$ we  shall write $Y_u \in 
\mathcal{X}(K/H)$ the associated fundamental vector field.
If we write $\pi \colon K \to K/H$ the 
canonical projection then we have
$\pi_* X_u = Y_{u}.$

\zati {\em The action $\Gamma \colon G \cap H \times H \to H$, defined 
by  $\Gamma(g,h) = g \cdot 
h$}.

\nt 
The action is tame since it extends to  the action $\Gamma \colon H \times H 
\to H$ 
defined by $\Gamma(g,h) = g \cdot 
h$.
 We write $\F_{H}$ the induced foliation, we have $\dim \F_{H} = \dim G \cap H$. 
For each $u \in \mathfrak{g} \cap \mathfrak{h}$ we  shall write $Z_u \in 
\mathcal{X}(H)$ the associated fundamental vector field.
Here $\mathfrak{h}$ is the Lie algebra of $H$.

\prg{\bf Twisted product.} This is the first geometrical tool we use for 
the study of $\mathcal{F}$.

 Take $K$ a connected compact abelian Lie group and $G$ a connected subgroup. We shall 
 suppose that $G$ is dense on $K$. Consider 
an orthogonal action  
$\Theta \colon H \times \R^n \to \R^n  $ where  $H$ is a compact subgroup of $K$ . 
The {\em twisted 
product} 
$K \times_{_H }\R^n$ is the quotient  of $K \times \R^n$ by the 
equivalence relation
$ (k\cdot  h^{-1},\Theta(h,
z)) \sim (k,z).
$
The element of the twisted product corresponding to $(k,z)$ is denoted 
by
$ <k,z>$. It is a 
manifold endowed with the action
$$
\Phi \colon G \times (K \times_{_H }\R^n) \TO (K \times_{_H }\R^n),
$$
defined by $\Phi (g, <k,z>) = <g \cdot k,z>.$ 
It is clearly a tame action
and we denoted by $\mathcal{F}_{tw}$ the induced 
conical foliation. 

The restriction $\Theta \colon G \cap H \times \R^n \to \R^n$  is 
also a tame action. The
induced conical foliation  is denoted by $\mathcal{F}_{\R^n}$;  a tamer 
group is  given by $H$. The couple 
$(\R^n,\mathcal{F}_{\R^n})$ is a {\em slice} of the twisted product.

The product $K \times \R^n$ is endowed with  the conical foliation 
$\F_{K} \times \F_{\R^n}$. The natural projection $R \colon K 
\times \R^n \to K \times_{H} \R^n$ gives the relations $  \F_{tw} = 
R_{*}(\F_{K} \times \F_{\R^n})$ and
\begin{equation}
    \label{strat}
\SFtw = \{ R ( K\times S) \ / \ S \in \SFRn 
\} = R(\{ K\} \times \SFRn).
\end{equation}

\prg {\bf Tubular neighborhood}. This is the second geometrical tool we use for 
the study of $\mathcal{F}$.

Consider a singular  stratum $S$ of $\SF$.  Since 
$\dim G(k\cdot x) = \dim G(x)$ for each 
$k \in K$ and each $x\in S$ then $S$ is  a $K$-invariant proper 
submanifold of $M$.
It possesses a $K$-invariant {\em tubular neighborhood } $
(T,\tau,S,\R^n)$ whose structural group is $O(n)$. We mean that $\tau$ 
is $K$-invariant and that there exists an atlas $\mathcal{B} $ such that 
for two charts $(U,\phii)$, $(V,\psi) \in \mathcal{B}$ and for a $k \in K$ 
the composition
$$
\psi \rondp k \rondp \phii^{-1} \colon (U \cap (k^{-1}   V))\times \R^n
 \TO  ((k  U )\cap    V) \times 
\R^n
$$
is of the form $ (x,u) \mapsto (k   \cdot x, A_{x,k}(u))$ with $A_{x,k} \in 
O(n)$.

Recall that there are the following smooth maps associated with this 
neighborhood:
\begin{itemize}
    \item[+] The {\em radius map} $\rho \colon T \to [0,1[$ defined fiberwise 
    from the assignation
$[x,t] \mapsto t$. Each $t\not= 0$ is a regular value of the $\rho$. 
The pre-image $\rho^{-1}(0)$ is $S$. This map is $K$-invariant, that is, 
$\rho(k \cdot z) = \rho(z)$.
     \item[+] The {\em contraction} $H \colon T \times [0,1] \to T$ 
     defined fiberwisely from $([x,t],r ) \mapsto  [x,rt]$. The restriction 
     $H_t \colon T \to T$ is an 
     embedding for each 
     $t\not= 0$   and $H_0 \equiv \tau$. We shall write $H(z,t) = t 
     \cdot z$. This map is $K$-invariant, that is, $ t \cdot (k \cdot z) = 
     k \cdot (t \cdot z)$.
\end{itemize}
These two maps are related by $\rho(t \cdot z) = t\rho(z)$. 

The hypersurface $D = \rho^{-1}(1/2)$ is the {\em 
     tube} of the tubular neighborhood. It is a $K$-invariant submanifold of $T$. 
     Notice that the map 
$$
\nabla \colon D \times [0,1[ \TO T,
$$
defined 
     by $\nabla(z,t) = (2t) \cdot z$ is a $K$-equivariant smooth map, 
     where $K$ acts trivially on the $[0,1[$-factor. Its restriction 
$\nabla \colon D \times ]0,1[ \TO T-S$ is a $K$-equivariant diffeomorphism
and its restriction $\nabla \colon D \times \{ 0 \} \equiv D \to S$ is 
$\tau$.

\smallskip

The foliation $\mathcal{F}$ induces on the fibers of $\tau$ a tame
foliation. Let us  see that. Write $G_S = G_{x_0}$ the isotropy 
subgroup of a point (and therefore, any point)  $x_0 \in S$. This group
acts effectively on $\tau^{-1}(x_0)$. The trace on $\tau^{-1}(x_0)$ of $\mathcal{F}$ is given by the 
action of $G_S$. Using a chart of the atlas $\mathcal{B}$ we identify 
$\tau^{-1}(x_0)$ with $\R^n$. The induced foliated manifold 
$(\R^n,\F_{\R^n})$ is the {\em slice} of the tube. 

The action of $G_S$ induces the 
orthogonal  action $\Theta \colon G_S\times \S^{n-1} \to \S^{n-1}$. This action  is 
an effective tame action. We shall write $\mathcal{G}_S$ the induced conical 
foliation. In 
fact, $(\S^{n-1},\mathcal{G}_S)$ is the link of $S$. The formula   \refp{1Dic} becomes 
\begin{equation}
\label{2Dic}
\dim \mathcal{F} = \dim \mathcal{F}_S + \dim \mathcal{G}_S.
\end{equation}
\medskip

The tubular neighborhood gives rise to a local blow up:

\bP
\label{comd}
Given a conical chart $(U,\phii) \in \mathcal{B}$ there exists a 
commutative diagram
\begin{equation}
\label{blo}
\begin{CD}
\R^{m-n}\times \S^{n-1}\times [0,1[ @>\phii'>> D\times [0,1[\\
    @V P VV @V\nabla VV\\
\R^{m-n}\times c\S^{n-1}@>\phii>> T,\\
    \end{CD}
    \end{equation}
where $\phii' \colon (\R^{m-n} \times \S^{n-1}\times [0,1[ 
,\mathcal{H} \times \mathcal{G} \times \mathcal{I})\TO (D\times 
[0,1[,\mathcal{F} \times \mathcal{I})$ is a foliated embedding.
\eP
\pro
The existence of the embedding $\phii'$ and the commutativity 
$\tau\rondp \phii' = \phii \rondp P$
is guaranteed by Lemma \ref{lema}. We also know that the restriction 
$\phii'  = \nabla_{_{min}}^{-1} \rondp \phii \rondp P \colon (\R^{m-n} \times \S^{n-1} \times ]0,1[ 
,\mathcal{H} \times \mathcal{G} \times \mathcal{I})\TO (D\times 
]0,1[,\mathcal{F} \times \mathcal{I})$ is a foliated embedding. Applying  
Lemma \ref{embed} we 
get that 
$\phii'  \colon (\R^{m-n} \times \S^{n-1} \times [0,1[ 
,\mathcal{H} \times \mathcal{G} \times \mathcal{I})\TO (D\times 
[0,1[,\mathcal{F} \times \mathcal{I})$ is a foliated embedding. 
\qed

\prg {\bf Molino's blow up}. 
The Molino' blow up \cite{Mo} of the foliation $\mathcal{F}$ produces a 
new foliation $\widehat{\mathcal{F}}$ of the same kind but of smaller depth. 
The main idea is to replace each point of closed strata by the fiber of a 
convenient tubular neighborhood. But in order to avoid the boundary that 
appears in
this procedure, we take the double.

 We suppose $\depth \SF > 0$.
Denote $S_{_{min}}$ the union of closed (minimal) strata
and choose $ 
T_{_{min}}$ a disjoint family of $K$-invariant tubular neighborhoods  of the 
 closed strata. The union of associated tubes is denoted by $D_{_{min}}$. 
Notice that the induced map $\nabla_{_{min}} \colon D_{_{min}} \times ]0,1[ \TO 
 T_{_{min}} -S_{_{min}}$ is a $K$-equivariant diffeomorphism.
 The {\em blow up} of $M$ is the manifold
$$
\widehat{M} = 
\left\{
\Big( D_{_{min}} \times ]-1,1[\Big) \coprod \Big( (M- S_{_{min}}) \times \{ 
-1,1\}\Big)
\right\} \Big/ \sim ,
$$
where $(z,t) \sim (\nabla_{_{min}} (z, |t|), t/|t|)$, and the map
$
\mathcal{L} \colon  \widehat{M} \TO M
$
defined by 
$$
\mathcal{L}(v) = 
\left\{
\begin{array}{ll}
\nabla_{_{min}} (z, |t|)& \hbox{if } v = (z,t) \in  D_{_{min}} \times 
]-1,1[ \\[,2cm]
z& \hbox{if }  v = (z,j) \in (M- S_{_{min}}) \times \{ 
-1,1\}.
\end{array}
\right.
$$
Notice that $\mathcal{L}$ is a continuous map whose restriction 
$
\mathcal{L} \colon \widehat{M} - \mathcal{L}^{-1}(S_{_{min}}) \to (M-S_{_{min}})
$
is a $K$-equivariant smooth trivial 2-covering.

Since the map $\nabla_{_{min}}$ is $K$-equivariant then $\Phi$ induces the action 
$\widehat{\Phi} \colon K \times \widehat{M} \TO \widehat{M}$ by 
saying that the blow-up $\mathcal{L}$ is $K$-equivariant. 
The open 
submanifolds $\mathcal{L}^{-1}(T_{_{min}})$  and 
$\mathcal{L}^{-1}(T_{_{min}} - S_{_{min}})$
are clearly $K$-diffeomorphic to 
$D_{_{min}} \times ]-1,1[$ and $D_{_{min}} \times (]-1,0[ \cup ]0,1[)$ respectively.
Notice that the depth of the restrictions of $\mathcal{F}$ to $D_{_{min}}$, 
$T_{_{min}} - S_{_{min}} $ and  $M_{_{min}} - S_{_{min}} $ is strictly smaller than $\depth \SF$.

The 
restriction 
$\widehat{\Phi} \colon G \times \widehat{M} \TO \widehat{M}$ 
is a tame action, whose tamer is just $\widehat{\Phi} \colon K \times \widehat{M} \TO \widehat{M}$.
The induced foliation is $\widehat{\mathcal{F}}$.
The associated foliations $\mathcal{F}$ and $\widehat{\mathcal{F}}$ are 
related by $\mathcal{L}$ which is a foliated map. Moreover, if $S$ is a not 
minimal stratum of $\SF$ then there exists an unique  stratum $S'$ of 
$\SFjat$
such that $\mathcal{L}^{-1}(S) \subset S'$. The family $\{ S'  \ / \ S \in \SF\}$ covers $\widehat{M}$ 
and verifies the relationship: $S_1 \prec S_2 \Leftrightarrow S'_1 \prec 
S'_2$. We conclude the important property
$$
\depth \SFjat < \depth \SF.
$$For any perversity $\per{p}$ on $M$ we define the perversity $\per{p}$ on 
$\widehat{M}$ by $\per{p}(S') = \per{p}(S)$.

\prg{\bf Orbit type stratification.}
Following the action $\Phi \colon K \times M \to M$, the points of $M$ 
are classified by this equivalence relation:
$$
x \sim y \Longleftrightarrow K_x = K_y.
$$
The induced partition $\SFi$ is the {\em orbit type stratification} of 
$M$ (see for example \cite{Br0}). The elements of this stratification are connected 
$K$-invariants 
submanifolds, called {\em ot-strata}. This stratification is finer than the 
stratification $\SFbar$ defined by the action $\Phi \colon K \times M \to M$,
but it verifies similar properties, in particular, $(\SFi,\prec)$ is a 
poset with finite depth.

Since an ot-stratum is a $K$-invariant submanifolds then it  possesses an  invariant 
neighborhood. A blow up can be constructed as in the previous 
framework: the {\em ot-blow-up}.  We write $\mathcal{N} \colon \widetilde{M} \to M$ the ot-blow-up, 
which is a $K$-equivariant continuous map relatively to an action 
$\widetilde\Phi \colon K \times \widetilde{M} \to \widetilde{M}$.
We have 
\begin{equation}
\label{ot}
\depth  \hbox{\sf S}_{\widetilde{\Phi }}
<
\depth \SFi.
\end{equation}

\section{ Basic Intersection cohomology}
The basic cohomology of a foliated space is the right cohomological tool 
to study the transverse structure foliation. A conical foliation has a  
transverse structure which reminds the stratified pseudomanifolds of 
\cite{GM}, and for these kind of singular spaces the most 
adapted cohomological tool is the intersection cohomology.
In this section we mix up the two ingredients and we introduce the basic 
intersection cohomology.

\smallskip

For the rest of this section, we fix $(M,\mathcal{F})$ a conical foliated 
manifold.

\smallskip

A differential form defined on the regular stratum may have a wild behavior 
relatively to the singular strata. But there are some of which a good 
contact with the singular part. These are the perverse forms, and from 
them we are going to construct the basic intersection cohomology.

There are several ways to define perverse forms: using a system of 
tubular neighborhoods (cf. \cite{Bry}), using  a global 
blow up (cf. \cite{S}) \ldots  ; in this work we introduce a more intrinsical way, 
using  the local blow ups we already have seen.

\smallskip

We are going to deal with 
differential forms on   products $(\hbox{manifold}) \times [0,1[^p$; 
these forms
are restrictions of differential forms defined on $(\hbox{manifold})\times 
]-1,1[^p$.

\prg{\bf Perverse forms.} 
A differential form defined on the regular stratum may have a wild behavior 
relatively to  singular strata. But there are some of them with a good 
contact with the singular part. These are the perverse forms, and from 
them we are going to construct the basic intersection cohomology.
Roughly speaking, perverse forms are differential forms defined on 
the regular stratum $R_\mathcal{F}$ which are extendable through local blow 
ups.

The differential complex  $\lau{\Pi}{*}{\mathcal{F}}{M \times 
[0,1[^p}$ 
of {\em perverse forms} of $M \times 
[0,1[^p$ is introduced by induction 
on  $\depth  \SF$. 
When this depth is  0 then 
$$\lau{\Pi}{*}{\mathcal{F}}{M \times [0,1[^p} = 
\hiru{\Omega}{*}{M 
\times [0,1[^p}.
$$

Consider now the 
generic case. 
A perverse form of $M \times 
[0,1[^p$ is  first of all a differential form $\om$  defined on  
$R_\mathcal{F} \times [0,1[^p \subset M \times [0,1[^p$ such that, for  a conical 
chart  $(U,\phii)$, there exists a perverse form 
$\om_\phibas$ of $\R^{m-n} \times 
\S^{n-1} \times [0,1[^{p+1}$  with
\begin{equation}
\label{PF}
\om_\phibas = (P_{\phibas}\times 
\hbox{identity}_{[0,1[^p})^*\omega \ \ \ \hbox{ on } 
\ \ \ \R^{m-n} \times R_{\mathcal G} \times ]0,1[ \times [0,1[^p.
\end{equation}
Notice that this condition makes sense since the restriction of the local 
blow up
$$
P_{\phibas} \colon
\underbrace{\R^{m-n} \times R_{\mathcal G} \times ]0,1[}_{\hbox{reg. 
stratum of } \R^{m-n}\times c\S^{n-1} }\TO 
\underbrace{ U \cap R_\mathcal{F}}_{\hbox{reg. stratum of } U} \
$$
is a 
diffeomorphism.

We  notice that 
\begin{equation}
    \label{pfi}
P_{\phibas}^* \colon \lau{\Pi}{*}{\mathcal{F}}{U } \TO 
\lau{\Pi}{*}{\mathcal{H}\times \mathcal{G} \times \mathcal{I}}{\R^{m-n} \times 
\S^{n-1} \times [0,1[} 
\end{equation}
is a  differential graded 
commutative algebra (dgca in short)  isomorphism.

\smallskip

The complex $\lau{\Pi}{*}{\mathcal{F}}{M }$ is a dgca since 
$(\om +\eta)_\phibas = \om_\phibas + \eta_\phibas$, 
$(\om \wedge\eta)_\phibas = \om_\phibas \wedge \eta_\phibas$ and
$(d\om)_\phibas = d\om_\phibas$.

\small

\prgg {\bf Remarks}.

\Zati The notion of perverse form depends on the foliation 
$\mathcal{F}$ through the stratification $\SF$.

\zati A local blow 
up produces a factor $[0,1[$. Further 
desingularisation would produce extra 
$[0,1[$ factors. This is the reason for introducing directly $[0,1[^p$.

\zati The perversion condition does not depend on the choice of the conical 
chart. In fact, if the condition \refp{PF} is satisfied for a conical chart 
then it is also satisfied for another conical  chart (cf. 
Proposition \ref{lema}).


\zati The local blow up of the cone $c\S^{n-1}$ is essentially  the product $\S^{n-1} \times 
[0,1[$. So, we have the natural identification $\lau{\Pi}{*}{\mathcal{G}}{\S^{n-1} \times [0,1[} = 
\lau{\Pi}{*}{c\mathcal{G}}{c\S^{n-1} }$ given by the assignation $\om \mapsto 
\om_{|R_{\mathcal{G} }\times ]0,1[}$.

\zati Given a tubular neighborhood $\tau \colon T \to S$ of a singular 
stratum (cf.2.3), the map $\nabla \colon D \times [0,1[ 
\to T$  gives an isomorphism between 
$ \lau{\Pi}{*}{\mathcal{F}}{T}$ and $\lau{\Pi}{*}{\mathcal{F} 
\times \mathcal{I}}{D 
\times [0,1[} $ (cf. Proposition \ref{comd} and Proposition \ref{lema}). 

\zati Let us illustrate this notion with an example. Consider the isometric action
$\Phi \colon \R \times \S^7 \to \S^7$ defined by $\Phi (t,(z_0,z_1,z_2,z_3) ) = 
(e^{a\pi i t}z_0,e^{b\pi i t}z_1,e^{c\pi i t}z_2,z_3)$ with $(a,b,c)\not= 
(0,0,0)$.  Recall 
that the induced foliation $\mathcal{F}$ is conical (cf. 1.2.1 (e)). We can see $\S^7$ as the 
join $\S^5 * \S^1$ where $\R$ acts freely on the first factor, inducing the foliation
$\mathcal{G}$, and  trivially on the second factor. There is just one singular stratum, 
namely $\S^1$, whose link is $(\S^5,\mathcal{G})$. This stratum has a ``global conical chart'' 
$(\S^1
\times c\S^5,\mathcal{I} \times c\mathcal{G})$ then
$$
\lau{\Pi}{*}{\mathcal{G}}{\S^7} = \hiru{\Om}{*}{\S^5\times \D^2}.
$$

\zati There are differential forms on $R_\mathcal{F}$ which are not 
perverse. Any differential form $\omega$ of $M$ is perverse, that is, 
$\hiru{\Om}{*}{M} \subset \lau{\Pi}{*}{\mathcal{F}}{M}.$
 In 
  fact we have 
  $\omega_\phibas =  P_{\phibas}^*\om$  on $\R^{n-m} \times \S^{n-1}\times [0,1[$
  for any conical chart $(U,\phii)$.

\begin{quote}
\pro For the first part it is sufficient to consider a smooth function on $R_\mathcal{F}$ going to infinity when approaching to the singular 
part. For the second part we proceed by induction on on  $\depth  \SF$. 
Consider $\om \in \lau{\Pi}{*}{\mathcal{F}}{M \times [0,1[^p}$.
For  a conical chart $(U,\phii)$,  we have that
$\om_\phibas = P_{\phibas}^*\om \in \hiru{\Om}{*}{\R^{m-n} \times \S^{n-1}
\times [0,1[^{p+1}}$. Here, we apply the induction hypothesis
\qed 
\end{quote}

\zati An open  foliated embedding $f \colon 
(M_1,\mathcal{F}_1) \to (M_2,\mathcal{F}_2)$ between two foliated conical manifolds
induces the dgca operator $f^* \colon \lau{\Pi}{*}{\mathcal{F}_2}{M_2} \to 
\lau{\Pi}{*}{\mathcal{F}_1}{M_1}$. 

\begin{quote}
\pro The embedding $f$ preserves conical charts  and therefore perverse forms.
\qed
\end{quote}

\prg{\bf Perverse degree.} 
The amount of transversallity of a perverse form $\om \in 
\lau{\Pi}{*}{\F}{M}$ 
with respect to a  
singular stratum $S$
is measured by the perverse degree $||\om||_S$. We define first the {\em local 
perverse degree} $||\om||_U$ for  a conical chart $(U,\phii)$ of a point 
of $S$.

Notice that a local blow up replaces $U \cap S$ by  $\R^{m-n} \times R_\mathcal{G}  \times \{ 0 \}$
and that the restriction 
\begin{equation}
\label{bup}
P_{\phibas} \colon \R^{m-n} \times R_\mathcal{G}  \times \{ 0 \}  \TO U \cap S
\end{equation}
is (isomorphic to) a  trivial bundle. Since the differential form $\om$ is  a perverse 
form, then the differential
form  $\phii^*\om$ (living on 
$\R^{m-n} \times R_\mathcal{G} \times ]0,1[  $) extends to the 
differential form $\om_\phibas$ 
(living on 
$\R^{m-n} \times R_\mathcal{G} \times [0,1[  $). 
Roughly speaking, the perverse degree $||\om||_U$  is the vertical degree of  
$\om_\phibas$ relatively to the added part, that is, the fibration \refp{bup}. 

More precisely, when 
$\om_\phibas =0$ on $\R^{m-n} \times 
R_\mathcal{G}  \times \{ 0 \} $ we put
$||\om||_U = -\infty$ 
and in the other cases

$$
||\om||_U = \min \left\{ k \in \N \ \Bigg/ \ 
\begin{array}{c}
\omega_\phibas (u_0, \ldots , u_k, 
\ldots)
\equiv 0  \\
\hbox{for each family } \{ u_0, \ldots , u_k \} \\
\hbox{of tangent vectors to the fibers of}  \\ 
P_{\phibas} \colon \R^{m-n} \times R_\mathcal{G}  \times \{ 0 \}  \TO U \cap S
\end{array}
\right\}.
$$
This number does not depend on the choice of the conical 
chart.
\begin{quote}
\pro
Take $(U,\phii_1)$ and 
$(U,\phii_2)$ two foliated charts .  
 From Lemma \ref{lema} we have  
$\om_{\phibas_1} = \Phi_{1,2}^*\om_{\phibas_2}$ and that 
the restriction $\Phi_{1,2} \colon \R^{m-n} \times R_\mathcal{G}  \times \{ 0 \}\to 
\R^{m-n} \times R_\mathcal{G}  \times \{ 0 \}$ is a 
diffeomorphism.
This implies that $||\om||^{\phibas_1}_U = ||\om||^{\phibas_2}_U$.
\qed
\end{quote}

Finally, we define the {\em 
perverse degree} 
$||\omega||_S$ by
\begin{equation}
||\omega||_S =  \sup \left\{ || \om ||_U \ / \  (U,\phii) \hbox{ is a 
conical chart of a point of } S \right\}.
\end{equation}

For two perverse forms $\omega$ and $\eta$ and a singular stratum $S$ 
we have:
\begin{equation}
\label{prod}
||\omega + \eta||_S \leq \max \big\{ ||\omega||_S , ||  \eta||_S\big\} \ , \  
||\omega \wedge \eta||_S \leq  ||\omega||_S +||  \eta||_S .
\end{equation}

For a perverse form $\om$ the perverse degree is smaller of the usual 
degree and verifies 
\begin{equation}
\label{cota}
||\om||_S \leq  \dim R_{\mathcal{G}_{S}} = \codim_M S -1.
\end{equation}

\prgg{\bf Remarks}.

\Zati The notion of perverse degree depends on the foliation 
$\mathcal{F}$ through the stratification $\SF$.

\zati 
The perverse degree of a  differential form of $M$ is $-\infty$ or 0 (cf. 
3.1.1 (g)). 

\zati Consider $\S^{n-1}$ endowed with a conical foliation without 
$0$-dimensional leaves and the disk 
$c\S^{n-1}$ endowed with the induced conical foliation. Then 
the perverse degree of a form $\om \in \lau{\Pi}{*}{c\mathcal{G}}{c\S^{n-1}} = 
\lau{\Pi}{*}{\mathcal{G} \times \mathcal{I}}{\S^{n-1}\times [0,1[}$ is just the degree of 
the restriction $\om_{|\S^{n-1} 
\times \{ 0 \}}$, 
where the degree of 0 is $ - \infty$.

\zati Given a  tubular neighborhood $\tau \colon T \to S$ (cf. 2.3) we know that we 
can identify  
$ \lau{\Pi}{*}{\mathcal{F}}{T}$ with  $\lau{\Pi}{*}{\mathcal{F} 
\times \mathcal{I}}{D \times [0,1[} $ through $\nabla$ (cf. 3.1.1 (e)).
Since $\nabla \colon D \times ]0,1[\to T-S$ is a $K$-diffeomorphism we 
have that the family of strata of $T$ is:
$
\{ \nabla(S' \times ]0,1[) \ / \ S' \hbox{ is a stratum of } D \} \cup \{ 
S \}.
$
So, for any $\om \in \lau{\Pi}{*}{\mathcal{F}}{T}$ we have
\begin{itemize}
\item $\izq\om\der_{\nabla(S' \times ]0,1[) } = ||\nabla^*\om||_{S' \times 
]0,1[}$.
\item $||\om||_S = \izq(\nabla^*\om)_{\big|(D \cap R_{\F}) \times \{ 0 \} }\der_\tau$, where 
$||\alpha||_\tau$ denotes the vertical degree of $\alpha \in 
\hiru{\Om}{*}{D \cap R_{\F}}$ relatively to the projection $\tau \colon (D \cap R_{\F}) \times \{ 
0 \} \equiv D \cap R_{\F} \to S$ (cf. Proposition \ref{comd}).
\end{itemize}

\zati Let us illustrate this notion with the example of 3.1.1 (f). The 
perverse forms are just the differential forms of $\S^5 \times \D^2$. The 
perverse degree $|| \ ||_{\S^1}$ is measured relatively to the trivial 
fibration $\S^5 \times \S^1 \to \S^1$. So, if we write 
$\theta_5$ 
a volume form of $\S^5$, $\theta_1$ 
a volume form of $\S^1$ and $(x,y)$ the coordinates of $\D^2$, we have
$$
||\theta_1||_{\S^1} = 0 \ \ ; \ \
||\theta_5||_{\S^1} = 5 \ \ ; \ \
|| \theta_5 \wedge (x dx + y dy)||_{\S^1} = -\infty.
$$

\zati A perverse form with $||\omega||_S \leq 0$ and $ || d\omega 
||_S \leq 0$ induces a 
differential form $\omega_S$ on $S$. When this happens for each  
stratum $S$ we conclude that $\omega \equiv \{ \omega_S\}$ is a Verona's {\em 
controlled form} (cf. \cite{V}).

\begin{quote}
\pro
Consider $(U_2,\phii_2)$ a conical chart of a point of a
singular stratum $S$. The conditions $||\omega||_S \leq 0$ and $ || d\omega 
||_S \leq 0$ give the existence of a form $\eta_{\phibas_{2}} \in 
\hiru{\Omega}{*}{U_2 \cap S}$ with $\omega_{\phibas_{2}} = P_{\phibas_{2}}^*  
\eta_{\phibas_{2}}$ 
on 
$\R^{m-n} \times R_\mathcal{G} \times\{ 0 \}$.
Let $(U_1,\phii_1)$ be another conical chart of a point of $S$ with $(U_1,\phii_1) \subset 
(U_2,\phii_2)$.
From the Proposition \ref{lema} we have 
$$
P_{\phibas_{1}}^*   \eta_{\phibas_1} = 
\omega_{\phibas_1} = 
\Phi_{1,2}^* \omega_{\phibas_2} = \Phi_{1,2}^* 
P_{\phibas_{2}}^* \eta_{\phibas_2}
=
P_{\phibas_{1}}^*  \eta_{\phibas_2}
$$ 
and therefore 
$ \eta_{\phibas_1}
=
 \eta_{\phibas_2}$ on $U_1 \cap S$.
This implies that 
the forms 
$\{ \eta_{\phibas_{i}} \}$ 
overlap and define a differential form $\omega_S \in \hiru{\Omega}{*}{S}$.
 \qed
\end{quote}

\zati An open foliated embedding $f \colon 
(M_1,\mathcal{F}_1) \to (M_2,\mathcal{F}_2)$ between two foliated conical manifolds
induces the dgca operator $f^* \colon \lau{\Pi}{*}{\mathcal{F}_2}{M_2} \to 
\lau{\Pi}{*}{\mathcal{F}_1}{M_1}$ which preserves the perverse degree 
(cf. Proposition \ref{lema}).

\prg {\bf Basic cohomology}.  The basic cohomology of a foliation 
$\mathcal{F}$  is an important  tool 
in the study its transversal structure and plays the r™le of the cohomology 
of the 
leaf space $\mf$, which can be a wild topological space. 

Consider $(M, \mathcal{F})$  a foliated 
manifold. 
A differential form $\om \in \hiru{\Om}{*}{M}$ is {\em basic} if 
$$
i_X \om = i_X d \om = 0,
$$
for each vector field $X$ on $M$ tangent to the foliation $\mathcal F$. For example, a 
function $f$ is basic iff $f$ is constant on the leaves. We shall denoted by 
$\hiru{\Om}{*}{\mf}$
the complex of basic forms. Since the sum and the product of basic 
forms are still basic forms, then the complex of basic forms is a 
dgca. Its cohomology $\hiru{H}{*}{\mf}$ is the {\em 
basic cohomology} 
of $(M,\mathcal{F})$.
We also use  the {\em relative basic cohomology} $\hiru{H}{*}{(M,F)/
\mathcal{F}}$, that is, the cohomology computed from the complex 
of basic forms vanishing 
on the saturated set $F$ .

When the foliation comes from a fibration $f \colon M \to B$ with 
connected fibers, then the leaf space $\mf$ is in 
fact the manifold $B$ 
and the basic cohomology $\hiru{H}{*}{\mf}$  is the cohomology of $B$.

\prgg{\bf Remarks}.

\Zati The basic cohomology does not use the stratification $\SF$. 

\zati Let us illustrate this notion with an example. Consider the isometric 
action
$\Psi \colon \R \times \S^5 \to \S^5$  defined by $\Psi (t,(z_0,z_1,z_2) ) = 
(e^{a\pi i t}z_0,e^{b\pi i t}z_1,e^{c\pi i t}z_2)$ with $(a,b,c)\not= (0,0,0)$. The induced 
foliation $\mathcal{G}$ is a regular one and we have that 
$\hiru{H}{i}{\S^5/\mathcal G}$ is

\begin{center}
\begin{tabular}{|c|c|c|c|c|}\hline
$i=0$ & $i=1$ & $i=2$ &$i=3$ &$i=4$\\ \hline
 $\R$&  0 & $\R \cdot [e]$ &  0 & $\R \cdot 
[e^2]$\\ \hline
\end{tabular}
\end{center}
Here, the cycle $e \in \hiru{\Om}{2}{\S^5/\mathcal G}$ is an {\em Euler 
form}.
Notice that this cohomology is finite dimensional and verifies the 
PoincarŽ Duality. These facts are always true for any regular isometric 
flow on a compact manifold (see \cite{EH} and \cite{ESH}).

Consider now the singular isometric flow defined in 3.1.1 (f).
Here we have that $\hiru{H}{i}{\S^7/\mathcal{F}}$ is
\begin{center}
\begin{tabular}{|c|c|c|c|c|c|c|}\hline
 $i=0$ & $i=1$ & $i=2$  & $i=3$ & $i=4$ & $i=5$ & $i=6$\\ \hline
  $\R$&  0 & 0 & 0 & $\R \cdot [\beta 
\wedge e] $ 
&0 & $\R 
\cdot [\beta \wedge e^2] $ \\ \hline
\end{tabular}
\end{center}
The cycle $\beta \in \hiru{\Om}{2}{\S^7/\mathcal F}$ is the form induced by 1 
from the 
double suspension $\S^7 \equiv \Sigma \Sigma \S^5$.

Notice that this cohomology is finite dimensional. This is true for any  isometric 
flow on a compact manifold (see \cite{W}, \cite{JI}).
On the other hand, the
PoincarŽ Duality is lost. We introduce in this work the basic intersection 
cohomology in order to recover this property.

\zati A smooth foliated map $f \colon  (M_1,\mathcal{F}_1) \to 
(M_2,\mathcal{F}_2)$  
induces  the dgca operator $f^* \colon \hiru{\Om}{*}{M_2  / 
\mathcal{F}_2}$ $\to
\hiru{\Om}{*}{M_1 / \mathcal{F}_1}$.
\medskip

\prg {\bf Basic intersection cohomology}.  
The stratification   $\SF$ induced by a conical foliation $\mathcal{F}$
is rich enough to support an intersection cohomology theory. 

Consider $(M, \mathcal{F})$  a 
conical foliated 
manifold.  A {\em perversity} is a map $\per{p} 
\colon \SiF \to \Z$, where $\SiF$ is the family of singular strata. There are two particular 
perversities we are going to use: 
\begin{itemize}
\item [-] the {\em constant perversity} $\per{i}$ defined by $\per{i}(S) 
=i$, where $i \in \Z$,
and 
\item[-] the {\em (basic) top 
perversity} $\per{t}$  defined by $\per{t}(S) = 
\codim_M \mathcal{F} - \codim_S \mathcal{F}_S -2 = \codim_M  S  - \dim 
\mathcal{G}_S -2 $. 
\end{itemize}
Any two perversities can be added.

Associated to an open foliated embedding $f \colon (M',\mathcal{F}') 
\to (M, \mathcal{F}) $ there exists a perversity on 
$(M',\mathcal{F}') $, still written $\per{p}$, defined by 
$\per{p}(S') = \per{p}(S)$ where $S' \in \SiF$ and 
$S \in \SiF$ with $f(S') \subset S$.

The basic intersection cohomology appears when one considers basic forms whose 
perverse degree is controlled by a given perversity.
We shall write
$$
\lau{\Om}{*}{\per{p}}{\mf}  = \left\{ \om \in 
\hiru{\Om}{*}{R_\mathcal{F} / \mathcal{F}_{R_\mathcal{F}}} \cap \lau{\Pi}{*}{\mathcal{F}}{M} \ / \ 
\max \left(||\om||_S ,||d\om||_S \right)\leq \per{p}(S) 
\ \ \forall S \in \SiF \right\}
$$
the complex of basic perverse forms whose perverse degree (and that of the their 
differential) is bounded by the perversity $\per{p}$. It is a differential 
complex, but it is not an algebra, in fact the wedge product acts in 
the following
way:
$$
\wedge \colon \lau{\Om}{i}{\per{p}}{\mf} \times 
\lau{\Om}{j}{\per{q}}{\mf} \TO 
\lau{\Om}{i+j}{\per{p} + \per{q}}{\mf}$$
(see  \refp{prod}).
The cohomology $\lau{\IH}{*}{\per{p}}{\mf}$ of this complex is the {\em 
basic intersection cohomology} of $(M,\mathcal{F})$, or BIC for short,
 relatively to the perversity $\per{p}$ (cf. \cite{SW2}).

\prgg {\bf Remarks}.

\Zati
The basic intersection cohomology coincides with the basic cohomology when 
the foliation $\mathcal{F}$ is regular, that is when $\depth \SF =0$. 
But it also generalizes the 
intersection cohomology of Goresky-MacPherson (cf. \cite{GM}) when the 
leaf space 
$B$ lies in the right category, that of stratified pseudomanifolds 
(cf. \cite{SW2}). 

\zati The perverse degree of a perverse form verifies \refp{cota}. But 
when this form is also a basic one we have
$$
||\om||_S \leq  \dim R_\mathcal{G} - \dim \mathcal{F}_{R_\mathcal{G} } = 
\codim_M \mathcal{F} - \codim_S \mathcal{F}_S-1 =(\per{t}+\per{1})(S)
$$
(cf. \refp{1Dic}).

\zati Si $\om \in \lau{\Om}{\ell-1}{\per{t}}{\mf}$ et 
$\phii \colon (\R^{m-n} \times c\S^{n-1},\mathcal{H} \times c\mathcal{G}) \TO 
(U,\mathcal{F})$ is a 
conical chart then 
$$
\om_{\phii} \equiv 0 \hbox{ \ \ on \ \ } \R^{m-n} \times 
R_{\mathcal{G}}\times \{ 0 \},
$$
where $\ell = \codim_{M}\F$.
\begin{quote}
    \pro Put $S$ the stratum of 
$\SF$ containing $\phii (\R^{m-n} \times \{\vartheta\})$. 
Notice that we have $\per{t}(S) = \ell- ((m-n )- \dim \mathcal{H})- 2 = 
n-\dim \mathcal{G} -2$. If the above assertion is not true then there exists 
$(x,z) \in 
\R^{m-n} \times R_{\mathcal{G}}$,
$
\{v_{1}, \ldots, v_{m-n-\dim \mathcal{H}}\} \subset T_{x}\R^{m-n}
$
and
$
\{w_{1}, \ldots, w_{n-1-\dim \mathcal{G}}\} \subset T_{z}\S^{n-1}
$
with 
$\om(x,z,0) (v_{1}, \ldots, v_{m-n-\dim \mathcal{H}}, w_{1}, \ldots, 
w_{n-1-\dim \mathcal{G}}) \neq 0$.
This   is not possible since the vectors
$
\{w_{1}, \ldots, 
w_{n-1-\dim \mathcal{G}}\}
$ are tangent to the fibers of $P_{\phibas}$ and
$
n- 1 -\dim \mathcal{G}
> \per{t}(S)$.
\qed
\end{quote}

\zati Consider a tubular neighborhood $\tau \colon T \to S$  (cf. 2.3)
and $\per{p}$ a perversity on $T$. From  3.2.1 
(d) we get that $\nabla^*$ stablishes an isomorphism between
 $
 \lau{\Om}{*}{\per{p}}{T/\mathcal{F}}
$
and
$$
\left\{ \omega \in \lau{\Om}{*}{\per{p}}{D\times  [0,1[
\Big/ \mathcal{F} \times \mathcal{I}}
\ / \ \izq \om_{\big|(D \cap R_{\F})}\der_\tau \leq \per{p}(S) \hbox{ and } 
\izq d\om_{\big|(D \cap R_{\F})}\der_\tau \leq \per{p}(S) \right\}.
$$ 
Here,  the perversity $\per{p}$  on 
$D \times [0,1[$ defined by 
$
\per{p}(S' \times [0,1[) = \per{p}(\nabla(S' \times ]0,1[))$.

\zati Let us illustrate this notion with the example of 3.1.1 (f) (see also 
3.3.1 (b)). A direct 
calculation gives that $\lau{H}{i}{\per{p}}{\S^7/\mathcal F}$ is

\begin{center}
\begin{tabular}{|c|c|c|c|c|c|c|c|}\hline
& $i=0$ & $i=1$ & $i=2$ &$i=3$ &$i=4$ & i=5 & i = 6\\ \hline
$\per{p} \leq  \per{-1}$ &  0 & 0 & $\R\cdot [\beta] $&  
0 & $\R \cdot [ \beta \wedge e]$ &  0 & $\R \cdot 
[\beta \wedge e^2]$\\ \hline
$\per{p} =\per{0}$ &  $\R$ & 0 &  0 &  
0 & $\R \cdot [\beta \wedge e]$ &  0 & $\R \cdot 
[\beta \wedge e^{2}]$\\ \hline
$\cdots$ & $\cdots$ & $\cdots$ & $\cdots $& $\cdots$ & $\cdots$ & 
$\cdots$  & $\cdots$\\ 
\hline 
$\per{p} =\per{t}$ &  $\R$ & 0 &  $\R \cdot [e]$ &  
0 & $\R \cdot [\beta \wedge e]$ &  0 &0\\ \hline
$\per{p} \geq \per{t}+\per{1} $ &  $\R$ & 0 &  $\R \cdot [e]$ &  
0 & $\R \cdot [e^2]$ &  0 &0\\ \hline
\end{tabular}
\end{center}

The first line is the  relative basic cohomology 
$\hiru{H}{*}{(\S^7,\S^1)/ 
\mathcal{F}}$,   the second line is the basic cohomology
$\hiru{H}{*}{\S^7 / 
\mathcal{F}}$ and the last line is the basic cohomology
$\hiru{H}{*}{(\S^7 - \S^1)/ 
\mathcal{F}}$. 
These cohomologies are finite dimensional and verify the following 
PoincarŽ Duality:
$$
\lau{\IH}{i}{\per{p}}{\mf} \cong \lau{\IH}{j}{\per{q}}{\mf},
$$
for $i+j =6$ and $\per{p} + \per{q} =\per{t}$.

We shall prove in this work that all these facts are always true for any 
abelian isometric 
flow on a compact manifold.

\zati Given an open  foliated embedding $f \colon  (M_1,\mathcal{F}_1) \to 
(M_2,\mathcal{F}_2)$ 
we have the  induced differential operator $f^* \colon \lau{\Om}{*}{\per{p}}{M_2 / 
\mathcal{F}_2} $ $\to 
\lau{\Om}{*}{\per{p}}{M_1/ \mathcal{F}_2}$ (cf.  3.1.1 (h) and   
3.3.1 (c)). If $f$ is a foliated diffeomorphism then $f^{*}$ is a 
differential
isomorphism.


\bigskip

We present now some of the technical tools used in this work. 
We fix $(M,\mathcal{F})$ a conical foliated manifold and $\per{p}$
a perversity.

\prg{\bf Local calculations}

The intersection basic cohomology, as the basic cohomology, 
is not easily computable. 
It becomes computable for singular foliations defined by an abelian isometric group as  will be 
observed in the next section. 
Nevertheless, the typical calculations for the BIC
are the 
classical ones.

\bP
\label{pro}
Let $(\R^k,\mathcal{H})$ be a simple foliation. Put $\per{p}$ the perversity 
defined on the conical foliated manifold $(\R^k \times M, \mathcal{H} 
\times \mathcal{F})$ by $\per{p}(\R^k \times S) = \per{p}(S)$.
The canonical projection $\pr \colon \R^k \times M \to M$ induces the isomorphism
$$
\lau{\IH}{*}{\per{p}}{\mf} \cong \lau{\IH}{*}{\per{p}}{\R^k \times M /  \mathcal{H} 
\times \mathcal{F}}.
$$
\eP

\pro
 Fix $a \in \R^k$ a basis point and put $\iota \colon M\rightarrow \R^k \times M$ the inclusion 
defined by $\iota (x) = (a,x)$. Notice that a conical chart $(U,\phii)$ on 
$M$ induces the conical chart $(\R^k \times U,\Ide_{\R^k} \times \phii)$ 
on $\R^k \times M$. Under these charts the projection $pr$ becomes the 
canonical projection $\R^k \times \R^{m-n} \times c\S^{n-1} \to \R^{m-n} \times 
c\S^{n-1}  $ defined by $(u,c,\zeta) \mapsto (v,\zeta)$. The inclusion $\iota$ 
becomes the inclusion $\R^{m-n} \times c\S^{n-1} \to \R^k \times \R^{m-n} \times 
c\S^{n-1} $ defined by $(v,\zeta) \mapsto (a,v,\zeta)$.
An inductive argument on the depth  gives that the operators 
 $\pr^* \colon \lau{\Om}{*}{\per{p}}{\mf} \to 
\lau{\Om}{*}{\per{p}}{ \R^k \times M/ \mathcal{H} \times \mathcal{F}}$ and 
$\iota^* \colon 
\lau{\Om}{*}{\per{p}}{\R^k  \times M/ \mathcal{H} \times \mathcal{F}} \to 
\lau{\Om}{*}{\per{p}}{\mf}$, are well-defined differential operators.

Since the composition $\iota^{*} \rondp \pr^*$ is the identity  then it suffices to 
prove that $\pr^{*} \rondp \iota^*$ is homotopic to the identity.

The foliated homotopy $k_1 \colon \R^k \times [0,1]
\rightarrow \R^k$ defined by $k_0(u,t) = tu$ induces the homotopy  $k_1 \colon \R^k \times M \times [0,1]
\rightarrow \R^k \times M $ defined by $k_1(u,x,t) = (k_0(u,t),x)$. 
This homotopy does not involve the $M$-factor, so it induces the 
morphism
$k_1^* \colon \lau{\Pi}{*}{\mathcal{H} \times \mathcal{F}}{\R^k\times M}
\to \lau{\Pi}{*}{\mathcal{H} \times \mathcal{F}}{\R^k\times M \times [0,1]}$ which preserves the 
perverse degree.
We also have $d \rondp k_1^* = k_1^* \rondp d$ and therefore the differential morphism
$$k_1^* \colon \lau{\Om}{*}{\per{p}}{\R^k\times M}
\to \lau{\Om}{*}{\per{p}}{\R^k\times M \times [0,1]}.$$

The  integration along the $[0,1]$-factor does not involves $M$. So, the
 operator 
$K \colon \lau{\Pi}{*}{\mathcal{H} \times \mathcal{F}}{\R^k\times M}
\to \lau{\Pi}{*-1}{\mathcal{H} \times \mathcal{F}}{\R^k\times M}$, given by $K \omega = 
\int_0^1  k_1^* \omega$, is well-defined and preserves the perverse degree.
On the other 
hand, it verifies the homotopy equality
$$
d \rondp K + K \rondp d = \pr^{*} \rondp \iota^* - \Ide.
$$

This implies  that $dK$ also preserves the perverse degree. We conclude that 
$$K \colon \lau{\Om}{*}{\per{p}}{\R^k  \times M/ \mathcal{H} \times \mathcal{F}} 
\to \lau{\Om}{* -1}{\per{p}}{\R^k  \times M/ \mathcal{H} \times 
\mathcal{F}} $$
is well-defined and that is a homotopy operator between $\pr^{*} \rondp 
\iota^*$ and the identity.
\qed

For the cone $(c\S^{n-1} ,c\mathcal{G})$ (cf. 1.1.1 (c)) we have:

\bP
\label{local5}
Let $\mathcal{G}$ be a conical foliation without 0-dimensional leaves on the sphere $\S^{n-1} $. A perversity 
$\per{p}$  on $c\S^{n-1} $ gives the perversity $\per{p}$ on $\S^{n-1} $ defined 
by $\per{p}(S) = \per{p}(S \times ]0,1[)$. The canonical projection 
$\chi\colon 
\S^{n-1}  \times 
]0,1[ \to \S^{n-1}  $ induces the isomorphism
$$
\lau{\IH}{i}{\per{p}}{c\S^{n-1}  / c\mathcal{G}} =
\left\{
\begin{array}{cl}
\lau{\IH}{i}{\per{p}}{\S^{n-1}  / \mathcal{G}} & \hbox{if } i \leq \per{p}(\{ 
\vartheta \}) \\
0 & \hbox{if } i > \per{p}(\{ \vartheta \}). 
\end{array}
\right.
$$
\eP

\pro
The statement about perversities is clear. From 3.2.1 (c) we have
$$
\lau{\Om}{j}{\per{p}}{c\S^{n-1}  / c\mathcal{G}} =
\left\{
\begin{array}{ll}
\lau{\Om}{j}{\per{p}}{\S^{n-1}  \times [0,1[ / \mathcal{G} \times \mathcal{I} }
& \hbox{if } j < \per{p}(\{ 
\vartheta \})\\
\lau{\Om}{j}{\per{p}}{\S^{n-1}  \times [0,1[/ \mathcal{G} \times \mathcal{I}} 
\cap d^{-1} \Ker \iota^* & \hbox{if 
} 
j = \per{p}(\{ \vartheta \})\\
\lau{\Om}{j}{\per{p}}{\S^{n-1}  \times [0,1[/ \mathcal{G} \times \mathcal{I}} 
\cap  \Ker \iota^* & \hbox{if } 
j > \per{p}(\{ \vartheta \})
\end{array}
\right.
$$
Here $\iota \colon \S^{n-1}  \to \S^{n-1}  \times [0,1[ $ is the inclusion defined 
by 
$\iota(\theta) = 
(\theta,0)$. Proceeding as in  the above Proposition we get the isomorphism
$
\lau{\IH}{i}{\per{p}}{c\S^{n-1}  / c\mathcal{G}} =
\left\{
\begin{array}{cl}
\lau{\IH}{i}{\per{p}}{\S^{n-1}  / \mathcal{G}} & \hbox{if } i \leq \per{p}(\{ 
\vartheta \}) \\
0 & \hbox{if } i > \per{p}(\{ \vartheta \}). 
\end{array}
\right.
$
\qed

\prg{\bf Mayer-Vietoris}.

 An open covering $\{ U , V \}$  of $M$ by saturated open 
subsets is a {\em basic covering}
when there exists a 
subordinated partition of the unity made up of  controlled basic  
functions. They 
may or may not exist.
For a such covering we have the Mayer-Vietoris short sequence
$$
0 \to \lau{\Om}{*}{\per{p}}{\mf} \to 
\lau{\Om}{*}{\per{p}}{U/\mathcal{F}} \oplus 
\lau{\Om}{*}{\per{p}}{V/\mathcal{F}}  \to
\lau{\Om}{*}{\per{p}}{(	U \cap V)/\mathcal{F}}  \to 0,
$$
where the map are defined by $\om \mapsto (\omega,\omega)$ and 
$(\alpha,\beta) \mapsto \alpha - \beta$. The third map is onto since the 
elements of the partition of the unity are controlled basic functions. Thus,  
the sequence
 is exact. This result is not 
longer true for 
more general coverings.   

The Mayer-Vietoris sequence allows to make computations when the conical 
foliated manifold is covered by finite suitable covering. The passage from 
the finite case to the general case may be done using an adapted version of 
the Bredon's trick of  \cite[page 289]{Br}:

\begin{BT}
    \label{Brt}
    Let $X$ be a paracompact topological space (resp. compact topological 
    space)
and  let $\{ U_\alpha\}$ be an open covering, closed for 
finite intersection. Suppose that $Q(U)$ is a statement about open subsets 
of $X$, satisfying the following three properties:
\begin{itemize}
\item[(BT1)] $Q(U_\alpha)$ is true for each $\alpha$;
\item[(BT2)] $Q(U)$, $Q(V)$ and $Q(U\cap V )$ $\Longrightarrow$ $Q(U\cup V 
)$, where $U$ and $V$ are open subsets of $X$; 
\item[(BT3)] $Q(U_i) \Longrightarrow Q\left({\displaystyle  \ \bigcup_i 
}U_i\right)$, 
where $\{ U_i\}$ is a disjoint family (resp. finite disjoint family) of open subsets of $X$.
\end{itemize}
Then $Q(X)$ is true.
\end{BT}
\prg {\bf Compact supports}. 

For the study of the PoincarŽ Duality we shall need the notion of 
cohomology with compact supports. We define the {\em support} of a perverse 
form $\om \in \lau{\Pi}{*}{\mathcal{F}}{M}$ as the closure (in $M$!) 
$$
\supp \om = \overline{\{ x \in M-\Sigma_\mathcal{F} \ / \ \om(x) \not= 0\}}.
$$
We have the relations $\supp (\om + \om') \subset \supp \om \cup \supp 
\om'$, 
$\supp (\om \wedge \om') \subset \supp \om \cap \supp 
\om'$ and $\supp d\om \subset \supp \om$. We denote 
$$
\lau{\Om}{*}{\per{p},c}{\mf} =
\left\{ \om \in \lau{\Om}{*}{\per{p}}{\mf} \ \Big/ \ \supp \om
 \hbox{ is compact} \right\}
$$
the complex of basic differential forms with compact support. It is a differential 
complex, but it is not an algebra, in fact the wedge product acts in this 
way:
$$
\wedge \colon \lau{\Om}{i}{\per{p},c}{\mf} \times 
\lau{\Om}{j}{\per{q}}{\mf} \TO 
\lau{\Om}{i+j}{\per{p} + \per{q},c}{\mf}$$
(see  \refp{prod}).
The cohomology $\lau{\IH}{*}{\per{p},c}{\mf}$ of this complex is the {\em 
basic intersection cohomology with compact support} of $(M,\mathcal{F})$, 
relatively to the perversity $\per{p}$. Of course, when $M$ is compact we 
have $\lau{\IH}{*}{\per{p},c}{\mf} = \lau{\IH}{*}{\per{p}}{\mf}$.

For $\Sigma_\mathcal{F} =\emptyset$, this notion generalizes the basic 
cohomology with compact supports. 

Given a basic covering $\{ U,V\}$ of $M$ we have the Mayer-Vietoris short sequence
$$
0 \to \lau{\Om}{*}{\per{p},c}{(U \cap V)/\mathcal{F}} \to 
\lau{\Om}{*}{\per{p},c}{U/\mathcal{F}_U} \oplus 
\lau{\Om}{*}{\per{p},c}{V/\mathcal{F}_V}  \to
\lau{\Om}{*}{\per{p},c}{\mf}  \to 0,
$$
where the map are defined by $\om \mapsto (\omega,\omega)$ and 
$(\alpha,\beta) \mapsto \alpha - \beta$. The third map is onto since the 
elements of the partition of the unity are controlled basic functions. Thus,  the 
sequence 	 
 is exact. 

We give now some local calculations. Given a simple foliation $(\R^k ,\mathcal{H}) 
\equiv (\R^a \times 
\R^b,\mathcal{J} \times \mathcal{I})$, where $\mathcal{J}$ is the one-leaf 
foliation of $\R^a$ and $\mathcal{I}$ the pointwise foliation of $\R^b$, we 
have $\lau{H}{*}{c}{\R^k/\mathcal{H}} = \R$ generated by 
$[f \ dx_1 \wedge \cdots \wedge dx_b]$ where $f$ is a {\em bump function}: 
 $f \in 
\mathcal{C}^\infty(\R^b)$ with ${\displaystyle \int_{\R^b} f =1}$ and  compact support.
\bP
\label{procomp}
Let $(\R^k,\mathcal{H})$ be a simple foliation. We have the isomorphism
$$
\lau{\IH}{*}{\per{p},c}{\mf} \cong 
\lau{\IH}{*+b}{\per{p},c}{\R^k \times M /  \mathcal{H} 
\times \mathcal{F}}
$$
given by  $[\beta] \mapsto [f \ dx_1 \wedge \cdots \wedge dx_b\wedge \beta]$,
where $f \in \mathcal{C}^\infty(\R^b)$ is bump function.
\eP
\pro
Notice  first $\lau{\Om}{*}{\per{p},c}{\R^k \times M /  \mathcal{H} 
\times \mathcal{F}} = \lau{\Om}{*}{\per{p},c}{\R^b\times M /  \mathcal{I} 
\times \mathcal{F}}$. It suffices to prove the case $b=1$.

Before executing the calculation let us introduce some notation. 
       Let $\beta$ be a differential form on $\hiru{\Om}{*}{R_{\F} \times M}$
       which does not 
       include the $dt$ factor. 
By $\Inte{-}{c}\beta (s)\wedge ds$ and
       $\Inte{c}{-}\beta  (s)\wedge ds$ we denote the forms on 
       $\hiru{\Om}{*}{R_{\F} \times M}$ 
       obtained from $\beta$ 
       by integration with respect to $s$:
       $\left(\Inte{-}{c}\beta (s)\wedge ds \right) 
       (x,t)(\vec{v}_{1}, 
       \ldots ,\vec{v}_{i}) = \Inte{t}{c}(\beta (x,s)(\vec{v}_{1}, 
       \ldots ,\vec{v}_{i}) ) \ ds $
      and   $\left( \Inte{c}{-}\beta (s)\wedge ds 
       \right)(x,t)(\vec{v}_{1}, 
       \ldots ,\vec{v}_{i}) )= \Inte{c}{t}(\beta 
       (x,s)(\vec{v}_{1}, 
       \ldots ,\vec{v}_{i}) ) \,  ds$
       where  $c \in ]a,b[$, $(x,t) \in R_{\F}\times ]a,b[$ and 
       $(\vec{v}_{1}, 
       \ldots ,\vec{v}_{i}) \in T_{(x,t)}( R_{\F} \times ]a,b[) $ .
      
A generic 
differential form  $\om$ of $\lau{\Om}{*}{\per{p},c}{\R\times M /  \mathcal{I} 
\times \mathcal{F}}$ is of the form
$$
\om = \alpha  + \beta \wedge dt,
$$
where $\alpha , \ \beta \in \lau{\Om}{*}{\per{p},c}{\R\times M /  \mathcal{I} 
\times \mathcal{F}}$ do not contain $dt$. Consider the differential 
operators
$$
\Delta \colon \lau{\Om}{*}{\per{p},c}{ \mf} \TO \lau{\Om}{*-1}{\per{p},c}{\R\times M /  \mathcal{I} 
\times \mathcal{F}} \ \ \hbox{ and } \ \ 
\nabla \colon \lau{\Om}{*}{\per{p},c}{\R\times M /  \mathcal{I} 
\times \mathcal{F}} \TO \lau{\Om}{*+1}{\per{p},c}{\mf},
$$
defined by 
$$
\Delta (\alpha  + \beta \wedge dt ) = \int_{-\infty}^\infty \beta (s) 
\wedge ds 
\ \ \hbox{ and } \ \ 
\nabla (\beta) = f \ dt  \wedge \beta.
$$
Notice that $\Delta \rondp \nabla =\Ide$ (up to a sign) 
which gives $\nabla^* \rondp 
\Delta^* = \Ide$ (up to a sign). We prove now $\Delta^* \rondp \nabla^* =\Ide$. Consider 
$[\om = \alpha + \beta \wedge dt] \in \lau{\IH}{*}{\per{p},c}{\R\times M /  \mathcal{I} 
\times \mathcal{F}} $. Define $\eta \in \lau{\Om}{*-1}{\per{p},c}{\R\times M /  \mathcal{I} 
\times \mathcal{F}} $ by
$$
\eta =
\left(
\int_{-\infty}^t f(s) \ ds
\right)
\left(
\int^{\infty}_t \beta(s) \wedge ds
\right)
-
\left(
\int^{\infty}_t f(s) \ ds
\right)
\left(
\int_{-\infty}^t \beta(s) \wedge ds
\right).
$$
A straightforward calculation gives 
$$
d \eta = f \ dt \wedge 
\int_{-\infty}^\infty \beta(s) \wedge ds
 - \beta \wedge dt - \alpha  =
\Delta^* \rondp \nabla^* (\om) - \om.
$$
This ends the proof.
\qed

For the computation of $\lau{\IH}{*}{\per{p}}{c\S^{n-1}  / c\mathcal{G}}$ we 
consider $g\in\mathcal{C}^\infty ([0,1[)$ with $g\equiv 1$ on $[0,1/4]$,
$g \equiv 0$ on $[3/4,1[$  
and ${\displaystyle \int_0^1}g=1$. 
\bP
\label{local6}
Let $\mathcal{G}$ be a conical foliation without 0-dimensional leaves on the sphere $\S^{n-1} $. 
We have the isomorphism 
$$
\lau{\IH}{i}{\per{p},c}{c\S^{n-1}  / c\mathcal{G}} =
\left\{
\begin{array}{cl}
0 & \hbox{if } i \leq \per{p}(\{ \vartheta \})+1\\
\lau{\IH}{i-1}{\per{p}}{\S^{n-1} / \mathcal{G}} & \hbox{if } i \geq\per{p}(\{ 
\vartheta \}) +2,
\end{array}
\right.
$$
given by $[g \ dt \wedge \om]\leftarrow [\om] $.
\eP
\pro
From 3.2.1 (c) we have
$$
\lau{\Om}{i}{\per{p},c}{c\S^{n-1}  / c\mathcal{G}} =
\left\{
\begin{array}{ll}
\lau{\Om}{i}{\per{p},c}{\S^{n-1}  \times [0,1[ / \mathcal{G} \times \mathcal{I} }
& \hbox{if } i< \per{p}(\{ 
\vartheta \})\\
\lau{\Om}{i}{\per{p},c}{\S^{n-1}  \times [0,1[/ \mathcal{G} \times \mathcal{I}} 
\cap d^{-1} \Ker \iota^* & \hbox{if 
} 
i= \per{p}(\{ \vartheta \})\\
\lau{\Om}{i}{\per{p},c}{\S^{n-1}  \times [0,1[/ \mathcal{G} \times \mathcal{I}} 
\cap  \Ker \iota^* & \hbox{if } 
i > \per{p}(\{ \vartheta \})
\end{array}
\right.
$$
Here $\iota \colon \S^{n-1}  \to \S^{n-1}  \times [0,1[ $ is the inclusion defined 
by 
$\iota(\theta) = 
(\theta,0)$. Consider a cycle $\om = \alpha + dt \wedge \beta \in 
\lau{\Om}{i}{\per{p},c}{c\S^{n-1}  / \mathcal{G}}$. Notice the equality
$$
\om   = d 
\left(\int_{-}^1 \beta (s) \wedge ds\right) ,
$$
with ${\displaystyle \left(\int_{-}^1 \beta (s) \wedge ds\right) } \in 
\lau{\Om}{i-1}{\per{p},c}{\S^{n-1}  \times [0,1[/ \mathcal{G} \times 
\mathcal{I}} .$ This gives $\lau{\IH}{i}{\per{p},c}{c\S^{n-1}  / c\mathcal{G}}  
=0 $ if $i \leq \per{p}(\vartheta)+1$. Now, it suffices to prove that the 
assignment $[\om] \mapsto [g \ dt \wedge \om]$ establishes an 
isomorphism between $\lau{\IH}{*-1}{\per{p}}{\S^{n-1}  / \mathcal{G}} $ and 
$\hiru{H}{*}{\lau{\Om}{\cdot}{\per{p},c}{\S^{n-1}  \times [0,1[/ 
\mathcal{G} \times \mathcal{I}} 
\cap  \Ker \iota^* } = \lau{\IH}{*}{\per{p},c}{(\S^{n-1}  \times [0,1[,\S^{n-1}  \times 
\{ 0 \})/\mathcal{G} \times \mathcal{I}}$.
The proof is exactly the same of the previous Proposition, replacing 
$-\infty$ by 0. \qed

\prg {\bf Twisted product}. 

We show in this section how to compute the BIC of the twisted product 
$K \times_{H} \R^n$  (cf. 2.2) in 
terms of the group $K$ and the slice $\R^n$. 

Notice 
first that a perversity $\per{p}$ on $(K \times \R^n,\F_{K} \times 
\F_{\R^n})$ (resp. $(K \times_{H} \R^n,\F_{tw})$) is determined by a perversity 
$\per{p}$ 
on $(\R^n,\F_{\R^n})$ by putting $\per{p}(K \times S) = 
\per{p}(S)$ (resp. $\per{p}(R(K \times S) )= 
\per{p}(S)$, see \refp{strat}. 

\bL
\label{pi}
The natural projection $R \colon K \times \R^n \to K \times_{H} 
\R^n $
induces the differential monomorphism
$$
R^* \colon \lau{\Om}{*}{\per{p}}{K \times_H \R^n/ \mathcal{F}_{tw}}
\TO 
\lau{\Om}{*}{\per{p}}{K \times \R^n/\mathcal{F}_K \times 
\mathcal{F}_{\R^n}},
$$
for any perversity $\per{p}$.
\eL
\pro
We proceed in several steps.

\Zati \underline{\em A  foliated atlas for 
 $\pi \colon K \to 
K/H$}.

Since $\pi$ is a $H$-principal bundle then it possesses an atlas 
$\mathcal{A}_{\#} =\left\{ 
\phii \colon \pi^{-1}(U) \TO U \times H \right\}$ made up with $H$-equivariant 
foliated charts. The $H$-equivariance means $\phii (h \cdot k) = (\pi (k), 
h\cdot h_0)$ if $\phii (k ) = (\pi (k), 
h_0)$. We study the foliation $\phii_*\mathcal{F}_{K}$.
This equivariance property gives $\phii_* X_u = (0,Z_u)$ for each $u 
\in \mathfrak{g} \cap \mathfrak{h}$. 
Thus, the trace of the foliation $\phii_*\mathcal{F}_{K}$ on the 
fibers of the canonical projection $\pr\colon U \times H \to U$ is 
$\mathcal{F}_H$.
On the other hand, since the map $\pi$ is a $G$-equivariant submersion then $\pi_* 
\mathcal{F}_{K} = \mathcal{F}_{K/H}$, which gives $\pr_*\phii_*\mathcal{F}_{K} = 
\mathcal{F}_{K/H}$. We conclude that $\phii_*\mathcal{F}_{K} \subset \mathcal{F}_{K/H} \times 
\mathcal{F}_H$. By dimension reasons we get $ \phii_*\mathcal{F}_{K} = \mathcal{F}_{K/H} \times 
\mathcal{F}_H$.
The atlas $\mathcal{A}_{\#}$ is an $H$-equivariant foliated atlas of 
$\pi$.

\zati \underline{\em A foliated atlas of $R \colon K \times \R^n \to K 
\times_{H} \R^n$}.

%
%
We claim that  $\mathcal{A}^{\#} = 
\{\overline{\phii} \colon \pi^{-1}(U) \times_H \R^n \TO U \times 
\R^n
\ / \ (U,\phii) \in \mathcal{A}_{\#} \}$ is a foliated atlas of 
$R$
where
the map 
$
\overline{\phii} 
$
is defined by $\overline{\phii}(<k,z>) =  
(\pi(k),(\Theta ( (\phii^{-1}(\pi(k),0))^{-1}  \cdot k ,  z)))$. This 
map  is a 
diffeomorphism whose inverse 
is $\overline{\phii}^{-1}(u,z) = <\phii^{-1}(u,0),z>$.  It verifies
$$
\overline{\phii}_*R_* (\mathcal{F}_K 
\times \mathcal{F}_{\R^n}) = 
\overline{\phii}_*R_* (\phii^{-1} \times \Ide_{\R^n})_*(
\mathcal{F}_{K/H} \times 
\mathcal{F}_H \times \mathcal{F}_{\R^n}).
$$
A straightforward calculation shows
$\overline{\phii} \rondp R \rondp (\phii^{-1} \times \Ide_{\R^n}) = 
(\Ide_U \times \Theta)$. 
Since $\mathcal{F}_H$ is defined by the action $\Gamma \colon G\cap H \times 
H \to H$ then $\Theta_* (\mathcal{F}_H \times \mathcal{F}_{\R^n}) = 
\mathcal{F}_{\R^n}$.
Finally we obtain
$
\overline{\phii}_*\mathcal{F}_{tw}=  \mathcal{F}_{K/H} \times \mathcal{F}_{\R^n}.
$


\zati \underline{\em Last step}.

 Given $(U,\phii) \in \mathcal{A}_{\#}$, we have the commutative diagram
$$
\begin{CD}
U \times H \times \R^n @>\phii>> K  \times \R^n\\
   @V P VV @V R VV\\
U \times \R^n @>\overline{\phii}>> K \times_H \R^n,\\
   \end{CD}
   $$
where $P(u,h,z) = (u, \Theta(h,z))$ and
$R^{-1}(\Ima \overline{\phii}) = \Ima \phii$.
%
 We claim that 
$$
R^* \colon \lau{\Om}{*}{\per{p}}{K \times_H \R^n/\F_{tw}}
\TO 
\lau{\Om}{*}{\per{p}}{K \times \R^n/\mathcal{F}_K \times \mathcal{F}_{\R^n}}
$$
is a well-defined morphism. Since it is a local question and we have $R^{-1}(\Ima \overline{\phii}) = \Ima \phii$, 
then it suffices to 
prove that the induced map $P^* \colon \lau{\Om}{*}{\per{p}}{U \times \R^n / \mathcal{F}_{K/H} \times 
\mathcal{F}_{\R^n}}
\TO 
\lau{\Om}{*}{\per{p}}{U \times H \times \R^n/ \mathcal{F}_{K/H} \times 
\mathcal{F}_H \times 
\mathcal{F}_{\R^n}}
$
is well-defined.
This comes from the fact that the map 
$$\nabla \colon (U \times H \times \R^n, \mathcal{F}_{K/H} \times 
\mathcal{F}_H \times 
\mathcal{F}_{\R^n}) \TO (U \times H \times \R^n, \mathcal{F}_{K/H} \times 
\mathcal{F}_H \times 
\mathcal{F}_{\R^n}),
$$
defined by $\nabla(u,h,z) = (u,h,\Theta(h,z))$ is a 
foliated diffeomorphism. Since $  \pr_0\rondp \nabla = R$, with $\pr_0 \colon 
U \times H \times \R^n \to U \times \R^n$ canonical projection then it suffices to 
apply 3.4.1 (f) and Proposition \ref{pro}.
The injectivity of $R^{*}$ comes from the fact that $R$ is a 
surjection.
\qed

\bL
\label{interior}
Let $\Phi \colon G \times M \to M$ be a tame action. Put $K$ a 
connected tamer group 
of $G$. We write $V_u \in \mathcal{X}(K)$ the fundamental vector field 
associated to an element $u$ of the Lie algebra of $K$. For any perversity 
$\per{p}$ the interior operator
$$
i_{V_u} \colon \left(\lau{\Om}{*}{\per{p}}{\mf} \right)^K
\TO \left( \lau{\Om}{*-1}{\per{p}}{\mf}\right)^K.
$$
is well-defined.
\eL
\pro
Since $K$ is connected then the vector field $V_u$ is 
$K$-invariant. So, the contraction operator $i_{V_u}$ preserves the 
$K$-invariants differential forms. Moreover, if $\om \in 
\left(\lau{\Om}{*}{\per{p}}{\mf} \right)^K$ and $X$ is a vector field on 
$M$ tangent to the foliation $\mathcal{F}$ then we have:
$$
i_Xi_{V_u}\om = -i_{V_u}i_X  \om= 0 \hbox{ and } 
i_Xdi_{V_u}\om = -i_Xi_{V_u} d\om = i_{V_u}i_X d\om  =0.
$$
 We end the proof if we show that 
$\om \in \left(\lau{\Om}{*}{\per{p}}{M \times [0,1[^p} \right)^K$
implies 
$i_{V_u}\om \in \lau{\Om}{*-1}{\per{p} }{M  \times [0,1[^p} $.
We proceed by induction on the depth of $ \SF$.

\medskip

{\bf 1. First step:} $\depth \SF =0$. The result is clear. 

\medskip

{\bf 2. Induction step.} Since the question is a local one we can consider 
$M= T$ a $K$-invariant tubular neighborhood of a singular stratum $S$.
Write $\wt{V_u} \in \mathcal{X}(D \times [0,1[)$ the fundamental vector 
field associated to $u$. This vector field is tangent to the boundary of 
$D \times [0,1[$,  write $U$ its restriction. On the other hand, since $\nabla \colon D \times [0,1[ \to T$ is a 
$K$-equivariant map, then $\nabla_*\wt{V_u} = V_u$.

We have, for each 
$\alpha \in \hiru{\Om}{*}{D}$, 
the relationship: $i_U \alpha (v_1, \ldots,v_j) \ne 0
\Rightarrow \alpha (v_1, \ldots,v_j) \ne 0$. This gives
the inequality $||i_U \alpha||_\tau \leq ||\alpha||_\tau $.

Denote $I = \Ide_{[0,1[^p}$.
We have 
$(\nabla \times I)^*\om \in\lau{\Om}{*}{\per{p}}{D \times [0,1[^{p+1}}$
with 
$\left|\left|\left( (\nabla \times I)^*\om\right)_{\big|D \times \{ 
0 \} \times
[0,1[^p}\right|\right|_\tau 
\leq \per{p}(S)$ and
$\left|\left|\left( (\nabla \times I)^*d\om\right)_{\big|D \times \{ 
0 \} \times [0,1[^p}\right|\right|_\tau 
\leq \per{p}(S)$ (cf. 3.4.1 (d)). By induction hypothesis we have
$$
(\nabla \times I)^*(i_{V_u}\om) =
 i_{\wt{V_u}}(\nabla \times I)^*\om \in \lau{\Om}{*}{\per{p}}{D \times [0,1[^{p+1}}.
$$
The result comes now from:
\begin{eqnarray*}
\izq\left( (\nabla \times I)^*(i_{V_u}\om) \right)_{\big|D\times \{ 
0 \} \times 
[0,1[^p}\der_\tau 
\! \!  \! \!  &=& \! \! 
\izq \left(i_{\wt{V_u}}(\nabla \times I)^*\om \right)_{\big|D \times \{ 
0 \} \times 
[0,1[^p}\der _\tau 
\! \! = \! \! 
\izq i_U \left( (\nabla \times I)^*\om )_{\big|D\times \{ 
0 \} \times 
[0,1[^p}\right)\der_\tau 
\\
\! \!  \! \! &\leq& \! \!  
\izq \left(\nabla \times I)^*\om \right)_{\big|D \times \{ 
0 \} \times [0,1[^p}\der_\tau 
\leq \per{p}(S),
\end{eqnarray*}
and
\begin{eqnarray*}
\izq\left( (\nabla \times I)^*(d(i_{V_u}\om) )\right)_{\big|D  \times \{ 
0 \} \times [0,1[^p}\der_\tau
&=& 
\izq-\left( (\nabla \times I)^*(i_{V_u}d\om )\right)_{\big|D \times \{ 
0 \} \times 
[0,1[^p}\der_\tau
\leq \per{p}(S)
\end{eqnarray*}
since $\om$ is $K$-invariant and 3.4.1 (d).
\qed

\prgg{\bf Fixing some notations about Lie algebras.}

Write $\mathfrak{k}$, $\mathfrak{g}$ and $\mathfrak{h}$ the Lie algebras 
of $K$, $G$ and $H$ respectively. Choose $\kappa$ an invariant riemannian 
metric on $K$, which exists by compactness.
Consider
$$\left\{ u_1, \ldots u_a,u_{a+1}, \ldots , u_{b}, u_{b+1},\ldots , 
u_c,u_{c+1},\ldots,u_f\right\}
$$
 an orthonormal basis of $\mathfrak{k}$
with $\left\{ u_1, \ldots u_b\right\}$ basis of $\mathfrak{g}$ and 
$\left\{ u_{a+1}, \ldots u_c\right\}$ basis of $\mathfrak{h}$.

For each index $i$ we write $X_i \in \mathcal{X}(K)$ the associated invariant vector 
field to $u_i$ (cf. 2.1.1 (a)).   
Let $\gamma_i \in \hiru{\Om}{1}{K}$  be the 
dual form of $X_i$, that is $\gamma_i = i_{X_i} \kappa$.  It is a cycle and it is
invariant by $K$, that is,  $k^* \gamma_i = \gamma_i$  for each $k\in K$. 
Since $K/H$ is an abelian Lie group then
$
\hiru{H}{*}{K/H}
= \bigwedge^* (\gamma_1,\ldots,\gamma_a,\gamma_{c+1},\ldots,\gamma_f).
$
The $\mathcal{F}_{K/H}$-basic differential forms on $
\bigwedge (\gamma_1,\ldots,\gamma_a,\gamma_{c+1},\ldots,\gamma_f)$ are exactly $
{\bigwedge}^*  
(\gamma_{c+1}, \ldots,\gamma_f)$.
This gives, 
\begin{equation}
\label{KH}
\hiru{H}{*}{K/H\Big/\mathcal{F}_{K/H}}=
{\bigwedge}^* (\gamma_{c+1},\ldots,\gamma_f) . 
\end{equation}

\bP
\label{isom}
Let $K \times_{H }\R^n$ be a twisted product. Let us suppose that the 
group $G$ is connected and dense in the  group $K$. Then
$$
\lau{\IH}{*}{\per{p}}{K \times_{_H }\R^n/\mathcal{F}_{{tw}}}
\cong
\hiru{H}{*}{K/H\Big/\mathcal{F}_{K/H}} 
\otimes
\left(\lau{\IH}{*}{\per{p}}{\R^n/\mathcal{F}_{\R^n}}\right)^H,
$$
for any perversity $\per{p}$.
\eP
\pro 
Since the operator
$
R^* \colon \lau{\Om}{*}{\per{p}}{K \times_H \R^n/\F_{tw}}
\TO 
\lau{\Om}{*}{\per{p}}{K \times \R^n/\mathcal{F}_K \times \mathcal{F}_{\R^n}}
$
is a monomorphism (cf. Lemma \ref{pi})
the it suffices to compute the cohomology of $\Ima R^{*}$.
We describe this complex in several steps. 

\medskip

$<i>$ \underline{Description of $\hiru{\Om}{*}{K \times R_{\mathcal{F}_{\R^n}}}$}. 

\smallskip

A differential form of $
\hiru{\Om}{*}{K \times R_{\mathcal{F}_{\R^n}}}$ 
is of the form
\begin{equation}
\label{form}
\eta + \sum_{1 \leq i_1 < \cdots < i_\ell \leq f}\gamma_{i_1} \wedge \dots 
\wedge \gamma_{i_\ell} \wedge \eta_{i_1, 
\ldots , i_\ell}, 
\end{equation}
where the forms $\eta, \eta_{i_1, 
\ldots , i_\ell}\in \hiru{\Om}{*}{K \times R_{\mathcal{F}_{\R^n}}}$ verify  
$i_{X_{i}} \eta = i_{X_{i}} \eta_{i_1, 
\ldots , i_\ell} = 0$ for all indices.

\medskip

$<ii>$  \underline{Description of 
$\lau{\Pi}{*}{\mathcal{F}_K \times \mathcal{F}_{\R^n}}{K \times \R^n}$}.

\smallskip

Since the foliation $\mathcal{F}_{K}$ is regular then we always can 
choose a conical chart of the form $(U_1 \times U_2, \phii_1 \times 
\phii_2)$ where $(U_1,\phii_1)$ is a foliated chart of $(K,
\mathcal{F}_K)$ and $(U_2,\phii_2)$ is  a conical chart of 
$(\R^n,\mathcal{F}_{\R^n})$. The local blow up of the chart 
$(U_1 \times U_2, \phii_1 \times 
\phii_2)$ is constructed from the second factor without modifying the 
first one. 
So, the differential forms $\gamma_i$ are always 
perverse forms and a
differential form 
$\om \in \lau{\Pi}{*}{\mathcal{F}_K \times \mathcal{F}_{\R^n}}{K \times \R^n}$ 
is of the form \refp{form}
where $\eta, \eta_{i_1, \ldots , i_\ell}\in 
\lau{\Pi}{*}{\mathcal{F}_K \times \mathcal{F}_{\R^n}}{K \times \R^n}$ verify  
$i_{X_{i}} \eta = i_{X_{i}} \eta_{i_1, 
\ldots , i_\ell} = 0$ for all indices.

\medskip

$<iii>$ \underline{Description of $\hiru{\Om}{*}{K \times R_{\mathcal{F}_{\R^n}}/
\mathcal{F}_K \times \mathcal{F}_{\R^n}}$}. 

\smallskip

Consider $\om \in \hiru{\Om}{*}{K \times R_{\mathcal{F}_{\R^n}}/
\mathcal{F}_K \times \mathcal{F}_{\R^n}}$.
Denote
$$\aleph\colon K \times (K \times \R^n) \TO 
(K \times \R^n) 
$$
the action defined by $\aleph(g, (k,u))  = (g\cdot k,u) $ 
and $\aleph_G \colon G \times (K \times \R^n) \TO 
(K \times \R^n) $ its restriction.
Since the fundamental vector fields of the action $\aleph_G$ are tangent to the foliation 
$\mathcal{F}_K \times \mathcal{F}_{\R^n}$, then $\om $ is $G$-invariant and, 
by density, also $K$-invariant.
So, this form $\om$ 
 lies on 
$
{\bigwedge}^*  
(\gamma_{1}, \ldots,\gamma_f)
\otimes \hiru{\Om}{*}{ R_{\mathcal{F}_{\R^n}}}.$ We need here the 
connectedness of $G$.

Since the $\mathcal{F}_{K}$-basic differential forms of 
${\bigwedge}^*  
(\gamma_{1}, \ldots,\gamma_f)$ are exactly
${\bigwedge}^*  
(\gamma_{b+1}, \ldots,\gamma_f)$, then we get the differential isomorphism
$$
\hiru{\Om}{*}{K \times R_{\mathcal{F}_{\R^n}}/
\mathcal{F}_K \times \mathcal{F}_{\R^n}}
\cong
{\bigwedge}^*  
(\gamma_{b+1}, \ldots,\gamma_f)
\otimes \hiru{\Om}{*}{ R_{\mathcal{F}_{\R^n}}/\mathcal{F}_{\R^n}}.
$$

\medskip

$<iv>$ \underline{Description of 
$\lau{\Om}{*}{\per{p}}{K \times \R^n / \mathcal{F}_K \times 
\mathcal{F}_{\R^n}}$}. 

\smallskip

From $<ii> $ and $<iii> $ it suffices to control the perverse degree of the 
forms 
$$
\eta + \sum_{b+1 \leq i_1 < \cdots < i_\ell \leq f}\gamma_{i_1} \wedge \dots 
\wedge \gamma_{i_\ell} \wedge \eta_{i_1, 
\ldots , i_\ell}  \in {\bigwedge}^*  
(\gamma_{b+1}, \ldots,\gamma_f)
\otimes \lau{\Pi}{*}{\mathcal{F}_{\R^n}}{ \R^n}.
$$
 Consider $S$ a stratum of $\SFRn$. From $||\gamma_i||_{K \times S} = 0$ 
and  $||\eta||_{K \times 
S} =  ||\eta||_{S} $,  we get $||\gamma_{i_1} \wedge \dots 
\gamma_{i_\ell} \wedge \eta_{i_1, 
\ldots , i_\ell}||_{K \times S}  =||\eta_{i_1, 
\ldots , i_\ell}||_S$. We conclude that  
$$
\lau{\Om}{*}{\per{p}}{K \times \R^n / \mathcal{F}_K \times 
\mathcal{F}_{\R^n}}
\cong
{\bigwedge}^*  
(\gamma_{b+1}, \ldots,\gamma_f)
\otimes 
\lau{\Om}{*}{\per{p}}{\R^n / \mathcal{F}_{\R^n}},
$$
as differential complexes.
\medskip

$<v>$ \underline{Description of 
$\Ima \left\{ 
R^* \colon \lau{\Om}{*}{\per{p}}{K \times_H \R^n/ \mathcal{F}_{tw}}
\TO 
\lau{\Om}{*}{\per{p}}{K \times \R^n/\mathcal{F}_K \times \mathcal{F}_{\R^n}}\right\}
$}.

\smallskip

We denoted by
$\{W_{a+1},\ldots , W_c\}$ the
fundamental vector fields of the action $\Theta \colon H \times \R^n \to 
\R^n$.
Consider now the action $\Upsilon\colon H \times (K \times \R^n) \to (K 
\times \R^n)$ defined by $\Upsilon(h,(k,z)) = ( k \cdot  h^{-1}, \Theta(h,z))$. 
This fundamental vector fields are $\{ (-X_{a+1}, W_{a+1}),
\ldots , (-X_{c},W_c)\}$.
Given $h \in H$ we write $\Upsilon_h \colon K \times \R^n  \to K \times 
\R^n $ the map defined by $\Upsilon_h 
(k,z) = \Upsilon(h,(k,z)) $.
Therefore
$$
\Ima 
R^*  =
\left \{ \omega \in  {\bigwedge}^*  
(\gamma_{b+1}, \ldots,\gamma_f)
\otimes 
\lau{\Om}{*}{\per{p}}{\R^n / \mathcal{F}_{\R^n}}.\ \Bigg/ \
\begin{array}{l}
\hbox{(i) } i_{X_i} \omega =  i_{W_i}\omega \hbox{ if  } 
a+1  \leq i \leq c \\[,3cm]
\hbox{(ii) } (\Upsilon_h)^*\omega = \omega \hbox{ for } h \in H  
\end{array}
\right\},
$$
(as differential complexes).
Recall that $h^*\gamma_i=\gamma_i$  for each $h \in H$ and each $i$.
The foliation $\mathcal{F}_{\R^n}$  is defined by the action of $G \cap H$ 
and
 then $i_{W_i}\omega = 0$ for $i \in \{ 
a+1,\ldots,b\}$. Since $\gamma_i(X_j) = \delta_{ij}$  then we 
get that $\Ima R^*$ is the differential complex
$$
\left \{ \omega \in {\bigwedge}^* 
(\gamma_{b+1}, \ldots,\gamma_{f})
\otimes \left( \lau{\Om}{*}{\per{p}}{\R^n/\mathcal{F}_{\R^n}} \right)^H\ \Big/ \
i_{X_i} \omega =  i_{W_i}\omega  
\hbox{ if  } 
b+1 \leq i\leq c
\right\} = 
$$
$$
{\bigwedge}^*  (\gamma_{c+1}, \ldots,\gamma_{f}) \otimes 
\underbrace{\left \{ \omega \in {\bigwedge}^*  
(\gamma_{b+1}, \ldots,\gamma_{c})
\otimes \left( \lau{\Om}{*}{\per{p}}{\R^n/\mathcal{F}_{\R^n}} \right)^H\ \Big/ \
 i_{X_i} \omega =  i_{W_i}\omega  
\hbox{ if  } 
b+1 \leq i\leq c
\right\} }_{A^*}.
$$
Consider the operator 
$
\Delta \colon \left( \lau{\Om}{*}{\per{p}}{\R^n/\mathcal{F}_{\R^n}} \right)^H
 \TO {\bigwedge}^*  
(\gamma_{b+1}, \ldots,\gamma_{c})
\otimes \left( \lau{\Om}{*}{\per{p}}{\R^n/\mathcal{F}_{\R^n}} \right)^H,
$
defined by
$$
\Delta (\beta) = 
 \beta + \sum_{b+1 \leq i_1 < \cdots < i_\ell \leq c} \ \gamma_{i_1} 
\wedge \dots 
\wedge \gamma_{i_\ell} \wedge 
(i_{W{_{i_\ell}} }\cdots i_{W{_{i_1}}}\beta).
$$
 A straightforward computation gives that the operator $\Delta$ is a 
 differential operator and that the restriction
$
\Delta \colon \left( \lau{\Om}{*}{\per{p}}{\R^n/\mathcal{F}_{\R^n}} \right)^H
 \TO A^*,
$
is an isomorphism. The inverse operator is 
$$\Delta^{-1}\left( \xi_0+ \sum_{b+1 \leq i_1 < \cdots < i_\ell 
\leq c}\gamma_{i_1} \wedge \dots 
\wedge \gamma_{i_\ell} \wedge \xi_{i_1, 
\ldots , i_\ell}\right) = \xi_0
$$
(cf. Lemma \ref{interior}).
We
conclude that the differential complex $\Ima R^*$ is isomorphic to
$$
{\bigwedge}^*  
(\gamma_{c+1}, \ldots ,\gamma_{f})
\otimes  \left( \lau{\Om}{*}{\per{p}}{\R^n/\mathcal{F}_{\R^n}} \right)^H
$$
 
$<vi>$ \underline{Last step}. 

\smallskip

Consider now the
the operator
$$
\hiru{H}{*}{\left( \lau{\Om}{*}{\per{p}}{\R^n/\mathcal{F}_{\R^n}} 
\right)^H} \hookrightarrow \left( 
\lau{\IH}{*}{\per{p}}{\R^n/\mathcal{F}_{\R^n}}\right)^H
$$
induced by the inclusion  
$
\iota \colon \left( \lau{\Om}{*}{\per{p}}{\R^n/\mathcal{F}_{\R^n}} 
\right)^H \hookrightarrow \lau{\Om}{*}{\per{p}}{\R^n/\mathcal{F}_{\R^n}} .
$
The usual arguments show that this operator is an isomorphism:
\begin{itemize}
\item[-] Monomorphism: $\omega = 
d\eta \Rightarrow \omega = d {\displaystyle \int_H} 
\Upsilon^*_h \eta \  dh $ 
\item[-] Epimorphism: $\omega  - \Upsilon^*_h 
\omega  = 
d\eta_h  \hbox{ for each } h \in H \Longrightarrow \omega  -  {\displaystyle  \int_H} 
\Upsilon^*_h \omega\  dh =
d {\displaystyle  \int_H} 
\eta _h \  dh $
\end{itemize}
We get 
$$
\lau{\IH}{*}{\per{p}}{K \times_{_H }\R^n/\mathcal{F}_{tw}}
\cong
\hiru{H}{*}{K/H\Big/\mathcal{F}_{K/H}} 
\otimes
\left(\lau{\IH}{*}{\per{p}}{\R^n/\mathcal{F}_{\R^n}}\right)^H,
$$
(cf. \refp{KH}).
\qed
\prgg {\bf Remark.}

The same procedure gives 
that 
the differential operator $\nabla = (\Ide_{{\bigwedge}^*
(\gamma_1,\ldots,\gamma_a,\gamma_{c+1}, \ldots,\gamma_{f})} 
\otimes \Delta^{-1})\circ R^*$ 
gives the  isomorphism
$$
\left(\lau{\Pi}{*}{\mathcal{F}_{tw}}{K \times_H \R^n}\right)^K \cong
{\bigwedge}^*
(\gamma_{1}, \ldots,\gamma_{a}, \gamma_{c+1},\ldots,\gamma_{f})
\otimes \left( \lau{\Pi}{*}{\mathcal{F}_{\R^n}}{\R^n} \right)^H
=
\hiru{H}{*}{K/H} \otimes \left( \lau{\Pi}{*}{\mathcal{F}_{\R^n}}{\R^n} \right)^H.
$$
\section{Cohomological properties of the BIC}

We prove in this section that the BIC of a compact foliated manifold 
$(M,\mathcal{F})$,
determined by an abelian isometric action $\Phi \colon G \times M \to M$,  is finite 
dimensional and verifies the PoincarŽ Duality. 

When the orbits 
of this action have the same dimension, that is when $\depth \SF =0$, 
then the foliation $\mathcal{F}$ is a (regular) riemannian foliation (cf. 
\cite{Mo}) and the BIC becomes the usual basic cohomology 
$\hiru{H}{*}{\mf}$. We already know from \cite{ESH} that $\hiru{H}{*}{\mf}$ is finite 
dimensional and verifies the PoincarŽ Duality. 
 For the generic case we are going to proceed by induction on 
the depth of $\SF$. 

\bigskip

But first of all we show how the BIC generalizes the usual basic cohomology.
The same situation appears for the intersection homology of a stratified 
pseudomanifold (cf. \cite{GM}).

\bp
\label{limit}
Let $(M,\mathcal{F})$ be a  foliated manifold determined by a tame 
action. Then
\begin{itemize}
    \item[i)] $\lau{\IH}{*}{\per{q}}{\mf} \cong \hiru{H}{*}{R_{\mathcal{F}}/ \mathcal{F}}$ 
if $ \per{q} > \per{t}$.
\item[ii)] $\lau{\IH}{*}{\per{0}}{\mf} \cong \hiru{H}{*}{\mf}$.
\item[iii)] $\lau{\IH}{*}{\per{p}}{\mf} \cong 
\hiru{H}{*}{(M,\Sigma_{\mathcal{F}})/ \mathcal{F}}$ 
if $\per{p} < 0$.
\end{itemize}
\ep

\pro The map $\omega 
\mapsto \omega$ gives the differential operator $I_M  
\colon  \lau{\Om}{*}{\per{q}}{\mf}
\to  \hiru{\Om}{*}{R_\mathcal{F}/\mathcal{F}}$.
The restriction map $\om \mapsto \om_{R_{\mathcal{F}}}$ defines the differential operators  $J_M  
\colon  
\hiru{\Om}{*}{\mf} 
\to \lau{\Om}{*}{\per{0}}{\mf}$  and $
K_M 
\colon 
\hiru{\Om}{*}{(M,\Sigma_{\mathcal{F}})/\mathcal{F}} \to 
\lau{\Om}{*}{\per{p}}{\mf}$ (cf. 3.2.1 (b)). 
We prove by induction on the depth of $\SF$ the following assertions.
\begin{itemize}
    \item[]$\mathfrak{A}_1(M,\mathcal{F}) =``I_M$ is a quasi-isomorphism". 
      \item[]$\mathfrak{A}_2(M,\mathcal{F}) =``J_M$ is a quasi-iso\-morphism". 
        \item[]$\mathfrak{A}_3(M,\mathcal{F}) =``K_M$ is a quasi-iso\-morphism". 
    \end{itemize}

\medskip

{\bf 1. First step:} $\depth \SF =0$. The singular part 
$\Sigma_\mathcal{F}$ is empty and therefore 
$\lau{\IH}{*}{\per{q}}{\mf}  = \lau{\IH}{*}{\per{0}}{\mf} 
= \lau{\IH}{*}{\per{p}}{\mf} =
\hiru{H}{*}{\mf}$ with $I_M =J_M =K_M = \Ide$.

\smallskip

{\bf 2. Induction step:} 
The family $\left\{M - 
 S_{_{min}},T_{_{min}}\right\}$ is a 
 $K$-invariant open covering of $M$ (cf. 2.4). Choose $\alpha \colon [0,1[ \to \R$ a 
smooth map with $\alpha \equiv 1$ on $[0,1/4[$  and $\alpha \equiv 0$ on 
$[3/4,1[$.  Write $f = \alpha \rondp 
 \rho_{_{min}}  \colon M \to \R$, which is a $K$-invariant map and therefore 
 $\mathcal{F}$-basic. Since $\supp f \subset T_{_{min}}$
and
$\supp (1-f) \subset M - 
 S_{_{min}}$ we conclude that the covering is a basic one. From 3.6 we have 
 a Mayer-Vietoris sequence and we get
$$
\mathfrak{A}_i(T_{_{min}}- S_{_{min}}, \mathcal{F}),
\mathfrak{A}_i(M - S_{_{min}},\mathcal{F}) \hbox{ and }
\mathfrak{A}_i(T_{_{min}}, \mathcal{F})
\Longrightarrow
\mathfrak{A}_i(M ,\mathcal{F}),
$$
for $i=1,2,3$. The induction hypothesis gives 
$\mathfrak{A}_i(T_{_{min}}- S_{_{min}}, \mathcal{F})$
and 
$\mathfrak{A}_i(M - S_{_{min}},\mathcal{F})$, it remains to prove 
$\mathfrak{A}_i(T_{_{min}}, \mathcal{F})$.

From  3.2.1 (d) we 
know that we can identify the perverse forms of $T_{_{min}}$ with the 
perverse forms of $D_{_{min}} \times [0,1[$. This identification sends 
basic forms to basic forms and preserves the 
perverse degrees relatively to any stratum different from those of $S_{_{min}}$. 
The perverse degree of a perverse form $\om $ of $D_{_{min}} \times [0,1[$
relatively to $S_{_{min}}$ becomes the vertical degree of 
the restriction $\om$ relatively to
$
\nabla_{_{min}} \equiv  \tau_{_{min}} \colon D_{_{min}} \times \{ 0 \} 
\equiv D_{_{min}} \TO S_{_{min}}$ (cf. 3.2.1 (d)). That is, 

$$
\begin{array}{ccl}
  \lau{\Om}{*}{\per{q}}{T_{_{min}} /\mathcal{F}} & \hbox{becomes} &
    \lau{\Om}{*}{\per{q}}{D_{_{min}}  \times 
[0,1[/\mathcal{F}\times \mathcal{I}}   \  \hbox{ (cf. 3.4.1 (b))} \\[,3cm]
  \lau{\Om}{*}{\per{0}}{T_{_{min}} /\mathcal{F}} & \hbox{becomes} &
  \{ \om \in \lau{\Om}{*}{\per{p}}{D_{_{min}}  \times 
[0,1[/\mathcal{F}\times \mathcal{I}}\ / \  \om_{|D _{_{min}} }= \tau_{_{min}}^{*} 
\eta \hbox{ with } \eta  \in \hiru{\Om}{*}{S_{_{min}}/\mathcal{F} } \} \\[,3cm]
 \lau{\Om}{*}{\per{p}}{T_{_{min}} /\mathcal{F}} & \hbox{becomes} &
    \{ \om \in \lau{\Om}{*}{\per{p}}{D_{_{min}}  \times 
[0,1[/\mathcal{F} \times \mathcal{I}}\ / \  \om_{|D_{_{min}}  }= 0\}. 
\end{array}
$$
 Proceeding as in Proposition \ref{pro} we prove that 

$$
\begin{array}{llc}
   \lau{\IH}{*}{\per{q}}{T_{_{min}} /\mathcal{F}}& \cong &  
    \lau{\IH}{*}{\per{q}}{D_{_{min}} /\mathcal{F}}   \\[,3cm]
    \lau{\IH}{*}{\per{0}}{T_{_{min}} /\mathcal{F}}& \cong & 
    \hiru{H}{*}{S_{_{min}} /\mathcal{F}}\\[,3cm]
  \lau{\IH}{*}{\per{p}}{T_{_{min}} /\mathcal{F}}& \cong & 0
\end{array}
$$
 Notice that $I_{T_{_{min}} }$ 
is induced by $ \pr\rond \nabla_{_{min}}^{-1} \colon T_{_{min}} 
-\Sigma_{\mathcal{F}}\to D_{_{min}} 
\times ]0,1[ \to D_{_{min}} $ and $J_{T}$ becomes $\tau^{*}$.

On the other hand, since $\nabla_{_{min}}  \colon  (D_{_{min}} \times 
]0,1[,\mathcal{F}\times \mathcal{I}) \to (T_{_{min}}  - \Sigma_{\mathcal{F}}, 
\mathcal{F}) $ is a foliated diffeomorphism, 
we get
$$
\hiru{H}{*}{(T_{_{min}} -\Sigma_{\mathcal{F}})/\mathcal{F}} \cong 
\hiru{H}{*}{D_{_{min}}  \times 
]0,1[/\mathcal{F} \times \mathcal{I}} \stackrel{\ref{pro}} {\cong}
\hiru{H}{*}{D_{_{min}} /\mathcal{F}}, 
$$
 where the isomorphism is induced by $\pr \rond \nabla^{-1}$. 
 The  induction hypothesis gives that $I_{T_{_{min}} }$ is a 
 quasi-isomorphism.
 
The inclusion $\iota_{_{min}} \colon S_{_{min}}  \to T_{_{min}} $ and the projection 
$\tau_{_{min}}  \colon 
T_{_{min}}  \to S_{_{min}} $ are foliated maps with $\tau_{_{min}}  \rondp \iota _{_{min}} = 
\Ide$, so they induce the operators 
$
\iota^* \colon \hiru{\Om}{*}{S_{_{min}} /\mathcal{F}} \to 
\hiru{\Om}{*}{T_{_{min}} /\mathcal{F}}
$
and 
$
\tau^* \colon \hiru{\Om}{*}{T_{_{min}} /\mathcal{F}} \to 
\hiru{\Om}{*}{S_{_{min}} /\mathcal{F}}
$
verifying $\iota_{_{min}} ^* \rondp 
\tau_{_{min}} ^* = \Ide$. The composition $\iota_{_{min}}  \rondp \tau_{_{min}} $ is 
homotopic to the identity by a foliated homotopy. This homotopy is 
just $H \colon T _{_{min}} \times [0,1] \to T_{_{min}} $ (cf. 2.4). So, the operator 
$\tau_{_{min}} ^*$ induces the isomorphism 
 $$
    \hiru{H}{*}{T_{_{min}} /\mathcal{F}} \cong  \hiru{H}{*}{S_{_{min}} 
    /\mathcal{F}}.
    $$ 
This proves that  $J_{T_{_{min}} }$ is a 
 quasi-isomorphism.

Let $\omega \in \hiru{\Om}{*}{(T_{_{min}} ,\Sigma_{\mathcal{F}})/\mathcal{F}}$ be a cycle. The 
above homotopy operator gives the relation $\om = d 
{\displaystyle \int_{0}^{1}}H^{*}\om$. Since the homotopy $H$ 
preserves the foliation $\mathcal{F}$  then $\int_{0}^{1}H^{*}\om \in \hiru{\Om}{*}{(T,\Sigma_{\mathcal{F}})/\mathcal{F}}$
and therefore
$$
\hiru{H}{*}{(T_{_{min}} ,\Sigma_{\mathcal{F}})/\mathcal{F}} = 0.
$$
This proves that  $K_{T_{_{min}} }$ is a 
 quasi-isomorphism.
\qed

\prg{\bf Finiteness.}

We prove that the BIC of the conical foliation $\F$ induced by an abelian 
isometric action 
 on a compact manifold is finite dimensional.
We proceed by induction on the depth of $\SF$. In order to decrease 
the depth we use the blow up 
of Molino. This will lead us to the twisted product through the invariant tubular neighborhood.

\bP
\label{tub}
Let $(M,\mathcal{F})$ be  a conical foliated manifold determined by a tame 
action.
Consider $(T,\tau,S,\R^n)$  a $K$-invariant tubular neighborhood of a compact 
singular stratum $S$.
If the BIC of the slice $(\R^n,\mathcal{F}_{\R^n})$ is finite dimensional then the BIC of 
the tube $(T,\mathcal{F})$ is also finite dimensional.
\eP

\pro
We can suppose that $G$ is
connected and we fix a connected tamer group $K$ where $G$ is dense.

We consider the orbit type stratification induced by the action $\Phi 
\colon K \times S \to S$ (cf. 2.5). We prove by induction on $\depth \SFi$ the 
following statement.
    \begin{center}
	$\mathfrak{A}(T,\mathcal{F}) =$
``The BIC
	$
	\lau{\IH}{*}{\per{p}}{T/\mathcal{F}}
	$
	is finite dimensional for each perversity $\per{p}$."
	\end{center}

Fix $\per{p}$ a perversity. Recall that any $K$-invariant 
submanifold of $M$ inherits naturally a perversity (cf. 3.4), written 
also $\per{p}$. We proceed in two steps.

\medskip

{\bf 1. First step:} $\depth \SFi =0$. The isotropy subgroup of any point
of $S$ is a compact subgroup $H \subset K$. The orbit space $S/K$ is 
a manifold and the natural projection 
$\pi \colon S \to S/K$ is a locally trivial bundle with fiber $K/H$.
Fix $\{ U_\alpha\}$ a good open covering of $S/K$ (cf. \cite{Br}), which is 
closed for finite intersections. For an open subset $U \subset S/K$ we 
consider the statement
\begin{center}
$Q(U) = $ ``The BIC
$\lau{\IH}{*}{\per{p}}{\tau^{-1}\pi^{_1}(U)/\mathcal{F}}$
 is finite dimensional''.
\end{center}
Notice that $Q(S) = \mathfrak{A}(T,\mathcal{F}) $. We get the result if we 
verify the three conditions of the Bredon's trick (cf. Proposition 
\ref{Brt}).

\begin{itemize}
\item[(BT1)] Since $U_\alpha$ is contractile we  can identify  
$\pi^{-1}(U_\alpha)$ with $U_\alpha \times K/H$. The group $K$ acts by $k_0 \cdot (x, 
k_1H) = (x,k_0 k_1 H)$. Fix $\{ x_0 \}$ a basis point of $U_\alpha$ and 
identify $\{ x_0 \} \times K/H$ with $K/H$. The contractibility of the open 
subset $U_\alpha$ 
gives a $K$-invariant  $O(n)$-isomorphism between $(T,\tau,S,\R^n)$ 
and 
$(U_\alpha \times 
\tau^{-1}(K/H), 
\Ide_{U_\alpha} \times \tau, U_\alpha \times K/H,\R^n)$. 
Notice that, identifying $\tau^{-1}(eH)$ with $\R^n$, the map $<k,u> 
\mapsto k\cdot u$ realizes  a $K$ diffeomorphism between $\tau^{-1}(K/H)$ 
and the twisted product $K \times_H \R^n$.

The contractibility of $U_\alpha$ gives:
$$
\lau{\IH}{*}{\per{p}}{\tau^{-1}\pi^{_1}(U)/\mathcal{F}}
\cong
\lau{\IH}{*}{\per{p}}{U_\alpha \times  \tau^{-1}(K/H) /\mathcal{I} \times \mathcal{F}}
\stackrel{\refp{pro}}{\cong} {\lau{\IH}{*}{\per{p}}{ \tau^{-1}(K/H) /\mathcal{F}}
}
\cong 
\lau{\IH}{*}{\per{p}}{ K \times_H \R^n /\mathcal{F}},
$$
which is finite dimensional following Proposition \ref{isom}. 

\item[(BT2)] The lifting of a partition of the unity subordinated to the 
covering $\{ 
U,V\}$ of  $S/K$ gives a controlled and basic partition of the unity subordinated 
to  the covering $\{ 
\tau^{-1}\pi^{_1}(U),\tau^{-1}\pi^{_1}(V)\}$ of $T$ (cf. 3.2.1 (b)). This 
covering is a basic one. Now we apply 3.6.

\item[(BT3)] Clear since 
$\lau{\IH}{*}{\per{p}}{\left(\tau^{-1}\pi^{_1}\left({\displaystyle  \ 
\bigcup_{i = 0}^m}U_i\right)\right)/\mathcal{F}}
=
{\displaystyle \bigoplus_{i = 0}^m \lau{\IH}{*}{\per{p}}{\tau^{-1}\pi^{_1}(U_i)/\mathcal{F}}}.
$
\end{itemize}
Notice that the compactness asked in the statement of the Proposition is used 
here by considering a finite covering $\{ U_i \ / \ 1 \leq i \leq m\}$. 
For a non finite covering the property (BT3) could be false.
\medskip

{\bf 2. Induction step}. 
Denote $S^{^{min}}$ the union of closed (minimal) strata of $\SFi$
and choose $ 
T^{^{min}}$ a disjoint family of $K$-invariant tubular neighborhoods  of the 
 closed strata. The projection map is written $\tau^{^{min}} \colon 
 T^{^{min}} \to S^{^{min}}$. This set is not empty since $\depth \SFi > 0$. 
The union of associated 
tubes is denoted by $D^{^{min}}$. It is a compact $K$-invariant submanifold 
verifying 
\begin{equation}
\label{dm}
\depth \hbox{\sf S}_{\Phi \colon K \times D^{^{min}} \to D^{^{min}}}
< \depth \SFi.
\end{equation}
The induced map $\nabla^{^{min}} \colon D^{^{min}} \times ]0,1[ \TO 
 T^{^{min}} -S^{^{min}}$ is a $K$-equivariant diffeomorphism, trivial 
 action on the $]0,1[$-factor. The radius 
 map $\rho^{^{min}} \colon T^{^{min}} \to [0,1[$ is a $K$-invariant 
 map. 

The family $\left\{\tau^{-1}(S - 
 S^{^{min}}),\tau^{-1}(T^{^{min}})\right\}$ is a 
 $K$-invariant open covering of $T$. Choose $\alpha \colon [0,1[ \to \R$ a 
smooth map with $\alpha \equiv 1$ on $[0,1/4[$  and $\alpha \equiv 0$ on 
$[3/4,1[$.  Write $f = \alpha \rondp 
 \rho^{^{min}} \rondp \tau \colon T \to \R$, which is a $K$-invariant map and therefore 
 $\mathcal{F}$-basic. Since $\supp f \subset \tau^{-1}(T^{^{min}})$
and
$\supp (1-f) \subset \tau^{-1}(S - 
 S^{^{min}})$ we conclude that the covering is a basic one. From 3.6 we 
 get an exact Mayer-Vietoris sequence
$$
0 \to
\lau{\Om}{*}{\per{p}}{\tau^{-1}(T^{^{min}} - 
S^{^{min}})/\mathcal{F}}
\to
\lau{\Om}{*}{\per{p}}{\tau^{-1}(S - 
S^{^{min}})/\mathcal{F}}
\oplus
\lau{\Om}{*}{\per{p}}{\tau^{-1}(T^{^{min}})/\mathcal{F}}
\to
\lau{\Om}{*}{\per{p}}{T/\mathcal{F}}
\to
0.
$$
The Five Lemma gives
$$
\mathfrak{A}(\tau^{-1}(T^{^{min}} - S^{^{min}}),\mathcal{F}),  \ \ 
\mathfrak{A}(\tau^{-1}(T ^{^{min}}),\mathcal{F})  \ \hbox{ and } \ 
\mathfrak{A}(\tau^{-1}(S - S^{^{min}}),\mathcal{F}),
\Longrightarrow 
\mathfrak{A}(T,\mathcal{F}).
$$

\nt We check now these  three conditions.

\medskip

{\bf (a)} $\mathfrak{A}(\tau^{-1}(T^{^{min}} - S^{^{min}})\ \mathcal{F})$.
The $K$-equivariant diffeomorphism  $\nabla^{^{min}}$ produces 
by pull-back the commutative diagram
$$
\begin{CD}
\tau^{-1}(D^{^{min}}) \times ]0,1[ @>\widetilde{\nabla^{^{min}}} >> 
\tau^{-1}(T^{^{min}}  -S^{^{min}})\\
    @V{\tau \times \Ide_{]0,1[} }VV  @V\tau VV\\
D^{^{min}} \times ]0,1[ @>\nabla^{^{min}} >> T^{^{min}}  - S^{^{min}}
    \end{CD}
   $$
where $\widetilde{\nabla^{^{min}}} $ is a $K$-diffeomorphism. So, we get 
$$
\lau{\IH}{*}{\per{p}}{\tau^{-1}(T^{^{min}} - 
S^{^{min}})/\mathcal{F}} \cong
\lau{\IH}{*}{\per{p}}{\tau^{-1}(D^{^{min}}) \times ]0,1[ /\mathcal{F} 
\times \mathcal{I}} 
\stackrel{\refp{pro}}{\cong} {\lau{\IH}{*}{\per{p}}{\tau^{-1}(D^{^{min}})  
/\mathcal{F}}},
$$
which is finite from the induction hypothesis. (see \refp{dm}).

\medskip

{\bf (b)} $\mathfrak{A}(\tau^{-1}(T^{^{min}}),\mathcal{F}) $. The idea is 
the following. We prove that the inclusion 
$\tau^{-1}(S^{^{min}}) \hookrightarrow \tau^{-1}(T^{^{min}})$ induces an 
isomorphism 
\begin{equation}
\label{iso}
\lau{\IH}{*}{\per{p}}{\tau^{-1}(T^{^{min}})/\mathcal{F}) }\cong
\lau{\IH}{*}{\per{p}}{\tau^{-1}(S^{^{min}})/\mathcal{F})}.
\end{equation}
Now, since the depth of 
$\hbox{\sf S}_{\Phi \colon K \times S^{^{min}} \to S^{^{min}}}$ is $0$, 
it suffices to apply the first step.

The 
contraction $H^{^{min}} \colon T^{^{min}} \times [0,1] \to T^{^{min}}$ is 
a $K$ invariant map with $H^{^{min}}_0 =  \iota \rondp \tau^{^{min}}$ and 
$H^{^{min}}_1 = \Ide_{\tau^{-1}(T^{^{min}})}$, where $\iota \colon 
S^{^{min}} \hookrightarrow T^{^{min}}$ is the natural inclusion. Notice that $\tau^{^{min}} \rondp 
\iota = \Ide_{\tau^{-1}(S^{^{min}})}$. The map $H^{^{min}}$ is 
locally the map
$H^{^{min}} \colon U \times \R^m \times [0,1] \to U \times \R^m$ defined 
by $H^{^{min}}(x,v,t) = (x,tv)$.  

Consider the induced commutative 
diagram by pull-back
$$
\begin{CD}
\tau^{-1}(T^{^{min}}) \times [0,1] @>\widetilde{H^{^{min}}} >> 
\tau^{-1}(T^{^{min}} )\\
    @V \tau \times \Ide_{[0,1]} VV  @V\tau VV\\
T^{^{min}} \times [0,1] @>H^{^{min}} >> T^{^{min}} .
    \end{CD}
    $$
Put $\widetilde{\tau^{^{min}}}$ and $\widetilde\iota$ the pull backs 
of $\tau^{^{min}}$ and $\iota$ respectively.  We have 
$\widetilde{H^{^{min}}}_{_1} =  \widetilde\iota \rondp \widetilde{\tau^{^{min}}}$ and 
$\Ide = \widetilde{H^{^{min}}}_{_1} = \widetilde{\tau^{^{min}}}\rondp 
\widetilde\iota .$ 

The operator $\widetilde{H^{^{min}}}$ is $K$-invariant and therefore is a 
foliated morphism: $\widetilde{H^{^{min}}}^*\mathcal{F} = \mathcal{F} 
\times \mathcal{I}$.
It is locally of the form
$
\widetilde{H^{^{min}}} \colon U \times \R^m \times \R^n \times [0,1] \to 
U \times \R^m \times \R^n$
with
$
\widetilde{H^{^{min}}} (x,v,w,t) = (x,tv,w).
$
Since the stratification induced by $\SF$ is 
$
\{ U \times \R^m \times S' \times [0,1] \ / \ S' \in \SFRn\}
$
then the perverse condition and the perverse degree are read on the 
$\R^n$-factor. So, the induced operator
$$
\widetilde{H^{^{min}}}^* \colon
\lau{\Om}{*}{\per{p}}{\tau^{-1}(T^{^{min}}) /\mathcal{F}}
\to
\lau{\Om}{*}{\per{p}}{\tau^{-1}(T^{^{min}}) \times [0,1]/\mathcal{F} \times \mathcal{I}}
$$
is well-defined. The integration along the $[0,1]$-factor does not 
involve $\R^n$. So, the integration operator 
$$
K \colon
\lau{\Om}{*}{\per{p}}{\tau^{-1}(T^{^{min}}) /\mathcal{F}}
\to
\lau{\Om}{*-1}{\per{p}}{\tau^{-1}(T^{^{min}}) /\mathcal{F} },
$$
given by $K(\omega	) = \int_0^1 \widetilde{H^{^{min}}}^*\omega$,
is well-defined. On the other hand, it verifies the homotopy equality:
$$
d \rondp K + K \rondp d =  \left(\widetilde{H^{^{min}}}_1\right)^*  - 
\left(\widetilde{H^{^{min}}}_0\right)^* = 
\left(\widetilde{H^{^{min}}}_1\right)^* - \Ide.
$$

This gives ${\widetilde{\tau^{^{min}}}}^* \rondp \widetilde\iota ^*= \Ide$.
Since $\widetilde\iota ^* \rondp {\widetilde{\tau^{^{min}}}}^* = \Ide$ the we 
get \refp{iso}.

\medskip

{\bf (c)} $\mathfrak{A}(\tau^{-1}(S - S^{^{min}}),\mathcal{F})$.
The idea is to construct a $K$-invariant tubular neighborhood 
$(E,\nu,\widetilde{S},\R^n)$ and a 
$K$-equivariant commutative diagram 
\begin{equation}
\label{diag}
\begin{CD}
E @>\mathcal{M} >> 
T\\
    @V\nu VV  @V\tau VV\\
\widetilde{S}@>\mathcal{N}>>  S
    \end{CD}
\end{equation}
verifying
\begin{itemize}
\item[i)] The map $\mathcal{N} \colon\widetilde{S} \to S$ is the ot-blow-up 
of $S$ relatively to 
the action $\Phi \colon K \times S \to S$ and 
$\depth {\sf S}_{\widetilde{\Phi}} < \depth \SFi$, where 
$\widetilde{\Phi}$ is the induced action of $K$ on $\widetilde{S}$.

\item[ii)] The restrictions
$
\mathcal{M} \colon 
\mathcal{M}^{-1}\tau^{-1}(S -S^{^{min}}) \to \tau^{-1}(S-S^{^{min}})
$
and 
$
\mathcal{M} \colon 
\mathcal{M}^{-1}\tau^{-1}(T^{^{min}} -S^{^{min}}) \to \tau^{-1}(T^{^{min}}-S^{^{min}})
$
are two  trivial 2-coverings. 

\item[iii)]  The invariant tubular neighborhood
$
(\mathcal{M}^{-1}\tau^{-1}(T^{^{min}}), \nu, \mathcal{N}^{-1}(T^{^{min}}),\R^n)$ is 
$K$-equivariantly diffeomorphic to 
$(\tau^{-1}(D^{^{min}}),\tau\times \Ide,
 (D^{^{min}})\times ]-1,1[,\R^n)$ 
\end{itemize}
The conditions ii) and iii) give
\begin{equation}
\label{doble}
\begin{array}{rcl}
\mathfrak{A}(\tau^{-1}(S - S^{^{min}}),\mathcal{F}) &\Longleftrightarrow&
\mathfrak{A}(\mathcal{M}^{-1}\tau^{-1}(S -S^{^{min}}) 
,\mathcal{E})\\
\mathfrak{A}(\tau^{-1}(T^{^{min}} - S^{^{min}}),\mathcal{F}) &\Longleftrightarrow&
\mathfrak{A}(\mathcal{M}^{-1}\tau^{-1}(T^{^{min}} -S^{^{min}}) 
,\mathcal{E}) \\
\mathfrak{A}(\tau^{-1}(D^{^{min}}) ,\mathcal{F}) &\Longleftrightarrow&
\mathfrak{A}(\mathcal{M}^{-1}\tau^{-1}(T^{^{min}}) 
,\mathcal{E}),
\end{array}
\end{equation}
where
 $\mathcal{E}$ is the 
 foliation induced by the action of $K$ on $E$. 

The family $\left\{\mathcal{M}^{-1}\tau^{-1}(S - 
 S^{^{min}}),\mathcal{M}^{-1}\tau^{-1}(T^{^{min}})\right\}$ is a 
 $K$-invariant open covering of $E$. Choose $\alpha \colon [0,1[ \to \R$ a 
smooth map with $\alpha \equiv 1$ on $[0,1/4[$  and $\alpha \equiv 0$ on 
$[3/4,1[$.  Denote $f = \alpha \rondp 
 \rho^{^{min}} \rondp \mathcal{N} \rondp \nu \colon 
 E \to \R$, which is a $K$-invariant map and therefore 
 $\mathcal{E}$-basic, where  $\mathcal{E}$ is the 
 foliation induced by the action of $K$ on $E$. Since $\supp f \subset 
\mathcal{M}^{-1}\tau^{-1}(T^{^{min}})$
and
$\supp (1-f) \subset \mathcal{M}^{-1}\tau^{-1}(S - 
 S^{^{min}})$ we conclude that the covering is a basic one. From 3.6 we 
 get an exact Mayer-Vietoris sequence
$$
0 \to
\lau{\Om}{*}{\per{p}}{\mathcal{M}^{-1}\tau^{-1}(T^{^{min}} - 
S^{^{min}})/\mathcal{E}}
\to
\lau{\Om}{*}{\per{p}}{\mathcal{M}^{-1}\tau^{-1}(S - 
S^{^{min}})/\mathcal{E}}
\oplus
\lau{\Om}{*}{\per{p}}{\mathcal{M}^{-1}\tau^{-1}(T^{^{min}})/
\mathcal{E}} \to 
$$
$$
\to
\lau{\Om}{*}{\per{p}}{E/\mathcal{E}}
\to
0.
$$
\nt We check now these  three conditions.
\medskip	
The Five Lemma and \refp{doble} give
$$
\mathfrak{A}(\tau^{-1}(T^{^{min}} - 
S^{^{min}}),\mathcal{F}),  \ \ 
\mathfrak{A}(\tau^{-1}(D^{^{min}}),
\mathcal{F})  \ \hbox{ and } \ 
\mathfrak{A}(E,\mathcal{E})
\Rightarrow
\mathfrak{A}(\tau^{-1}(S - S^{^{min}}),\mathcal{F}) .
$$

\nt We check now these  three conditions.

\medskip

\begin{itemize}
\item[\bf (c1)] $\mathfrak{A}(\tau^{-1}(T^{^{min}} - 
S^{^{min}}),\mathcal{F})$. It is the condition {\bf (a)}.

\item[\bf (c2)] $\mathfrak{A}(\tau^{-1}(D^{^{min}}),
\mathcal{F}) $. By induction hypothesis since we have \refp{dm}.

\item[\bf (c3)] $\mathfrak{A}(E,\mathcal{E})$. By 
induction hypothesis since we have \refp{ot} by i).

\end{itemize}

\bigskip

It remains to construct \refp{diag}. Consider the manifold
$$
\widetilde{S} = 
\left\{
\Big( D^{^{min}} \times ]-1,1[\Big) \coprod \Big( (S- S^{^{min}}) \times \{ 
-1,1\}\Big)
\right\} \Big/ \sim ,
$$
where $(z,t) \sim (\nabla^{^{min}} (z, |t|), t/|t|)$, and the map
$
\mathcal{N} \colon  \widetilde{S} \TO S
$
defined by 
$$
\mathcal{N} (v) = 
\left\{
\begin{array}{ll}
\nabla^{^{min}} (z, |t|)& \hbox{if } v = (z,t) \in  D^{^{min}} \times 
]-1,1[ \\[,2cm]
z& \hbox{if }  v = (z,j) \in (S- S^{^{min}}) \times \{ 
-1,1\}.
\end{array}
\right.
$$
This is the ot-blow-up 
$\mathcal{N}\colon \widetilde{S} \to S$ induced by the action $\Phi 
\colon K \times S \to S$.

Consider the manifold
$$
E = 
\left\{
\Big( \tau^{-1}(D^{^{min}}) \times ]-1,1[\Big) \coprod \Big( 
\tau^{-1}(S - S^{^{min}}) \times \{ 
-1,1\}\Big)
\right\} \Big/ \sim ,
$$
where $(z,t) \sim (\widetilde{\nabla^{^{min}} }(z, |t|), t/|t|)$ the map
$
\mathcal{M} \colon  E \TO T
$
defined by 
$$
\mathcal{M} (v) = 
\left\{
\begin{array}{ll}
\widetilde{\nabla^{^{min}} }(z, |t|)& \hbox{if } v = (z,t) \in  \tau^{-1}(D^{^{min}}) \times 
]-1,1[ \\[,2cm]
z& \hbox{if }  v = (z,j) \in \tau^{-1}(S - S^{^{min}})  \times \{ 
-1,1\}
\end{array}
\right.
$$
and the map $\nu \colon E \to \widetilde{S}$ 
defined by 
$$
\nu (v) = 
\left\{
\begin{array}{ll}
(\tau(z),t)& \hbox{if } v = (z,t) \in  \tau^{-1}(D^{^{min}}) \times 
]-1,1[ \\[,2cm]
(\tau(z),j)& \hbox{if }  v = (z,j) \in \tau^{-1}(S - S^{^{min}})  \times \{ 
-1,1\}
\end{array}
\right. 
$$
Since   $\nabla^{^{min}}$ and   $\widetilde{\nabla^{^{min}} }$ are $K$-equivariant 
embeddings
then $\widetilde{S}$  and $E$ are $K$-manifolds. The maps 
$\mathcal{N}$ and $\mathcal{M}$ are $K$-equivariant continuous 
maps. Since the map $\tau$ is $K$-equivariant then the map $\nu$ is a $K$-equivariant map. 
The diagram \refp{diag} is clearly commutative.

We have that $(E,\nu,\widetilde{S}, \R^n)$ is a 
tubular neighborhood since
$(\tau^{-1}(S - S^{^{min}}),\tau,S - S^{^{min}}, \R^n)$
and
$(\tau^{-1}(D^{^{min}}) \times ]-1,1[ , \tau \times \Ide_{]-1,1[},D^{^{min}}\times ]-1,1[ , \R^n)$
are compatible tubular neighborhoods. It remains to verify the properties 
i)-iii).
\begin{itemize}
\item[i)] By construction.
\item[ii)] By construction $\mathcal{M}^{-1}\tau^{-1}(S -S^{^{min}}) = 
\tau^{-1}(S -S^{^{min}}) \times \{-1,1\}$ and $\mathcal{M}$ is the 
projection on the first factor.
\item[iii)] By construction $\mathcal{N}^{-1}(T^{^{min}}) = D^{^{min}} 
\times ]-1,1[$,
$\mathcal{M}^{-1}\tau^{-1}(T^{^{min}}) = \tau^{-1}(D^{^{min}} ) \times 
]-1,1[$ and $\nu$ becomes $\tau \times \Ide_{]-1,1[}$.
\end{itemize}
This ends the proof. \qed

The first main result of this section is the following 
\bT
\label{T}
The BIC of the foliation determined by an isometric action of an abelian 
Lie group on a compact 
manifold is finite dimensional.
\eT

\pro
Given a conical foliated manifold $(N\mathcal{N})$ we consider the 
statement 
\begin{center}
	$\mathfrak{A}(N,\mathcal{N}) =$
``The BIC
	$
	\lau{\IH}{*}{\per{p}}{N/\mathcal{N}}
	$
	is finite dimensional, for any perversity $\per{p}$."
	\end{center}

Consider $\Phi \colon G \times M \to M$ an isometric action of an abelian Lie group $G$
on a compact manifold. This equivalent to say that the action is tame. 
We denote by  
$\mathcal{F}$ the induced conical foliation. Let us suppose that $G$ is 
connected and 
dense in the (connected) tamer group $K$. 
We prove $\mathfrak{A}(M,\mathcal{F}) $ by induction on $\depth \SF$.
\medskip

{\bf 1. First step:} $\depth \SF =0$. The foliation $\mathcal{F}$ is a 
(regular) riemannian foliation (cf. \cite{Mo}) and the BIC is just the basic 
cohomology (cf. 3.4.1 (a)) The result 
comes directly from \cite{ESH}.

\smallskip

{\bf 2. Induction step:} 
The family $\left\{M - 
 S_{_{min}},T_{_{min}}\right\}$ is a 
 $K$-invariant basic open covering of $M$ (cf. proof of Proposition 
 \ref{limit}). From 3.6 we 
 get an exact Mayer-Vietoris sequence
$$
0 \to
\lau{\Om}{*}{\per{p}}{T_{_{min}} - 
S_{_{min}}/\mathcal{F}}
\to
\lau{\Om}{*}{\per{p}}{M - 
S_{_{min}}/\mathcal{F}}
\oplus
\lau{\Om}{*}{\per{p}}{T_{_{min}}/\mathcal{F}}
\to
\lau{\Om}{*}{\per{p}}{\mf}
\to
0.
$$

The Five Lemma gives
$$
\mathfrak{A}(T_{_{min}} - S_{_{min}},\mathcal{F}),  \ \ 
\mathfrak{A}(T _{_{min}},\mathcal{F}) \ \hbox{ and } \ 
\mathfrak{A}(M - S_{_{min}},\mathcal{F}),
\Longrightarrow 
\mathfrak{A}(M,\mathcal{F}).
$$

\nt We check now these  three conditions.

\medskip

{\bf (a)} $\mathfrak{A}(T^{^{min}} - S^{^{min}}, \mathcal{F})$.
Since $\nabla^{^{min}}$ is a $K$-equivariant diffeomorphism    we have
$$
\mathfrak{A}(T^{^{min}} - S^{^{min}}, \mathcal{F})
\Longleftrightarrow
\mathfrak{A}(D^{^{min}} \times ]-1,1[,\mathcal{F} \times \mathcal{I})
\stackrel{\refp{pro}}{\Longleftrightarrow}
\mathfrak{A}(D^{^{min}}, \mathcal{F} ),
$$
which is true since the depth of $(D^{^{min}},\mathcal{F})$ is strictly 
smaller than $\depth \SF$ (cf. 2.5).

\medskip

{\bf (b)} $\mathfrak{A}(T^{^{min}} , \mathcal{F})$. If we prove 
$\mathfrak{A}(\R^n , \mathcal{F}_{\R^n})$, for the slice of a tubular 
neighborhood $(T,S,\tau,\R^n)$, 
then it suffices to apply the
Proposition \ref{tub}. Recall that $\mathcal{F}_{\R^n}$ is defined by an 
orthogonal action $\Theta \colon G \cap H \times \R^n \to \R^n$ such 
 that $\Theta \colon G\cap H \times \S^{n-1} \to 
\S^{n-1}$ is a tame action without fixed points  defining 
$\mathcal{G}$ and verifying 
  $(\R^{n}, \mathcal{F}_{\R^{n}}) =(c\S^{n-1}, c\mathcal{G}_{S}).
$
We have
$
\lau{H}{*}{\per{p}}{\R^n/\mathcal{F}_{\R^n}} 
\stackrel{\ref{local5}}{\cong}
\lau{H}{\leq \per{p}(\vartheta)}{\per{p}}{ \S^{n-1} /\mathcal{G}_{S}} ,
$
the tronqued cohomology, 
which is finite dimensional since 
$
\depth {\sf S}_{\mathcal{G}_{S}} < \depth \SF
$
implies $\mathfrak{A}(\S^{n-1},\mathcal{G}_{S})$.

\medskip

{\bf (c)} $\mathfrak{A}(M - S^{^{min}}, \mathcal{F})$.
The family $\left\{\mathcal{L}^{-1}(M - 
 S_{_{min}}),\mathcal{L}^{-1}(T_{_{min}})\right\}$ is a 
 $K$-invariant open covering of $M$. Choose $\alpha \colon [0,1[ \to \R$ a 
smooth map with $\alpha \equiv 1$ on $[0,1/4[$  and $\alpha \equiv 0$ on 
$[3/4,1[$.  Write $f = \alpha \rondp 
 \rho_{_{min}}  \rondp \mathcal{L}\colon \widehat{M }\to \R$, 
which is a $K$-invariant map and therefore 
 $\widehat{\mathcal{F}}$-basic. Since $\supp f \subset 
 \mathcal{L}^{-1}(T_{_{min}})$
and
$\supp (1-f) \subset \mathcal{L}^{-1}(M - 
 S_{_{min}})$ we conclude that the covering is a basic one. From 3.6 we 
 get an exact Mayer-Vietoris sequence
$$
0 \to
\lau{\Om}{*}{\per{p}}{\mathcal{L}^{-1}(T_{_{min}} - 
S_{_{min}})/\widehat{\mathcal{F}}}
\to
\lau{\Om}{*}{\per{p}}{\mathcal{L}^{-1}(M - 
S_{_{min}})/\widehat{\mathcal{F}}}
\oplus
\lau{\Om}{*}{\per{p}}{\mathcal{L}^{-1}(T_{_{min}})/\widehat{\mathcal{F}}}
\to
\lau{\Om}{*}{\per{p}}{\widehat{M}/\widehat{\mathcal{F}}}
\to
0.
$$
Recall that, following 2.5, we have that $
\mathcal{L} \colon \mathcal{L}^{-1}(M - 
S_{_{min}})\to (M-S_{_{min}})
$
is a $K$-equivariant smooth trivial 2-covering. So, 
$\mathfrak{A}(\mathcal{L}^{-1}(M - 
S_{_{min}}),\widehat{\mathcal{F}}) \Longleftrightarrow 
\mathfrak{A}(M - 
S_{_{min}},\mathcal{F})$. Now, the Five Lemma gives
$$
\mathfrak{A}(\mathcal{L}^{-1}(T_{_{min}} - 
S_{_{min}}),\widehat{\mathcal{F}}),  \ \ 
\mathfrak{A}(\mathcal{L}^{-1}(T_{_{min}}),\widehat{\mathcal{F}}) \ \hbox{ and } \ 
\mathfrak{A}(\widehat{M},\widehat{\mathcal{F}}),
\Longrightarrow 
\mathfrak{A}(M - 
S_{_{min}},\mathcal{F}).
$$

\nt We check now these  three conditions.

\medskip

\begin{itemize}
\item[]
\begin{itemize}
\item[{\bf (c1)}] $\mathfrak{A}(\mathcal{L}^{-1}(T_{_{min}} - 
S_{_{min}}),\widehat{\mathcal{F}})$.
Since $\mathcal{L}^{-1}(T_{_{min}} - 
S_{_{min}})$ is $K$-diffeomorphic to two copies of $T_{_{min}} - 
S_{_{min}}$  (cf. 2.4)) then    we have
$$
\mathfrak{A}(\mathcal{L}^{-1}(T_{_{min}} - 
S_{_{min}}),\widehat{\mathcal{F}})
\Longleftrightarrow
\mathfrak{A}(T^{^{min}} - S^{^{min}}, \mathcal{F}).
$$
Now we apply (a).

\item[{\bf (c2)}] $\mathfrak{A}(\mathcal{L}^{-1}(T_{_{min}}),\widehat{\mathcal{F}}) $. 
From 2.4 we know that $\mathcal{L}^{-1}(T_{_{min}})$ is $K$-diffeomorphic 
to $D_{_{min}} \times ]-1,1[$. Now we proceed as in (a).

\item[{\bf (c3)}] $\mathfrak{A}(\widehat{M},\widehat{\mathcal{F}})$.
 Because $\depth {\sf S}_{\widehat{\mathcal{F}} }< \SF$ (cf.  2.4). \qed
\end{itemize}
\end{itemize}

\prg {\bf PoincarŽ Duality}.
 We prove in this section that the BIC of a conical 
 foliation  $\F$ defined on an oriented manifold 
$M$ and
determined by a tame action,
verifies the PoincarŽ Duality:
\begin{equation}
\label{pairing}
\lau{\IH}{*}{\per{p}}{\mf} \cong \lau{\IH}{\ell -*}{\per{q},c}{\mf}.
\end{equation}
Here  $\ell$ (or $\ell_M$) is the codimension of the foliation 
$\mathcal{F}$. 
The  two  perversities $\per{p}$ and 
$\per{q}$ are complementary, that is, $\per{p} + \per{q} = 
\per{t}$.

The proof follows the path of 4.2, but firstly we define the morphism $P_M$
giving \refp{pairing}, it depends on the notion of a tangent volume form.

\prgg {\bf Tangent volume form}. For the definition of the pairing 
$P_{M}$ we need a volume form tangent to the leaves of $\F$.

Consider $\Phi \colon G \times M \to M$ a tame action defining $\F$. We 
can choose the group $G$ connected and the action $\Phi$ effective; so, $b = \dim G = \dim 
\F$.
We also fix
$\{u_1,\ldots u_b\}$ a basis of the 
Lie algebra $\mathfrak{g}$ of $G$. The associated fundamental vector fields on $M$ are 
denoted by $\{ V_1, \ldots,V_b\}$. 

A {\em tangent volume form} of $(M,\mathcal{F})$ is a $G$-invariant differential form 
$\eta \in \lau{\Pi}{b}{\mathcal{F}}{M}$ verifying

\begin{equation}
    \label{eta}
\eta(V_1,\ldots,V_b) =1.
\end{equation}


\nt Notice that $d\eta(V_1,\ldots,V_b,-) =0$.

\bL
\label{etatang}
Consider $K$ a tamer group of $G$. There exists a $K$-invariant tangent volume form
$\eta$
of $(M,\mathcal{F})$ 
verifying the following properties:

\Zati  For each $\om \in \lau{\Om}{\ell}{\per{t}}{\mf}$ the 
product $\om \wedge \eta$ does not depend on the tangent volume form $\eta$.

\zati For each $\om \in \lau{\Om}{\ell-1}{\per{t}}{\mf}$ the 
product $\om \wedge d\eta$ is 0.

\zati For each $\om \in \lau{\Om}{\ell}{\per{t},c}{\mf}$ the 
integral $\displaystyle{\int_{R_{\mathcal{F}}} \omega \wedge \eta }$ is finite.

\zati For each $\om \in \lau{\Om}{\ell-1}{\per{t},c}{\mf}$ the 
integral $\displaystyle{\int_{R_{\mathcal{F}}} d(\omega \wedge \eta) }$ 
is 0.
\eL
\pro
We proceed in several steps.

\medskip

\nt {\em \underline{First step}. Existence.}

\smallskip

We prove the following statement by induction 
on $\depth \SF$:
\begin{quote}
    ``There exists a $K$-invariant differential form
    $\eta \in \lau{\Pi}{b}{\F \times \mathcal{I}}{M\times 
	[0,1[^p}$ verifying:
    \begin{equation}
	\label{etap}
	\eta((V_1,0),\ldots,(V_b,0)) =1,
	\end{equation}
	for each $p \in \N$''.
	\end{quote}
	The existence is proven by taking $p=0$.
	
	\smallskip
	
 When $\depth \SF =0$ then we define $\eta_{0}$ on 
the orbits of $\Phi$  by \refp{etap} and we extend to a differential form of
$\hiru{\Om}{b}{M \times [0,1[^p}$. The differential form 
$\eta = \Int{K}{}k^{*}\eta_{0}$ is $K$-invariant, 
	 lives in $\lau{\Pi}{b}{\F \times \mathcal{I}}{M\times 
     [0,1[^p}$ and verifies \refp{etap} since each $k$ is a 
     $K$-equivariant diffeomorphism.

Consider now the 
case $\depth \SF >0$. 
By induction hypothesis there exists  a $K$-invariant differential 
form $\eta_{0} 
\in \lau{\Pi}{b}{\wh\F}{\wh{M}\times 
    [0,1[^p}$ verifying \refp{etap}. Associated to 
the Molino's blow up we have the $
K$-equivariant imbedding
$
 S \colon (M- 
S_{_{min}}) \to \mathcal{L}^{-1}(M- 
S_{_{min}}),
$
defined by
$
\sigma(z) = (z,1).
    $
    The differential form $\eta = (\sigma\times \hbox{identity}_{[0,1[^p})^{*}\eta_{0}$ 
    belongs to $\hiru{\Om}{b}{R_{\F}\times 
    [0,1[^p}$. It is $K$-invariant and verifies \refp{etap} since 
    $\sigma$ is a $K$-equivariant imbedding. It remains to prove that
    $\eta \in \lau{\Pi}{b}{\F \times \mathcal{I}}{M \times 
    [0,1[^p}$, which is a local property.
So, we  can consider that $M$ is a tubular 
 neighborhood  $T$ of a singular stratum of $\SF$ and  prove 
 $(\nabla\times \hbox{identity}_{[0,1[^p})^{*} \eta \in 
\lau{\Pi}{b}{\F \times \mathcal{I}}{D \times [0,1[^{p+1}} $ (cf. 3.1.1 (e)).
This is the case since  $ \sigma \rondp \nabla \colon D \times ]0,1[ \to D \times ]-1,1[$ is just the 
 inclusion and $\eta_{0} \in \lau{\Pi}{b}{\F \times \mathcal{I}}{D\times 
    [0,1[^{p+1}}$.
\medskip

\nt {\em \underline{Second step}. The condition (a).}

\smallskip

Let $\eta'$ be another tangent volume form associated to 
$\mathcal{F}$ through $\Phi$ and $\{u_{1}, \ldots , u_{b}\}$. 
 By degree reasons it suffices to prove the equality
 $
 i_{V_{1}} \cdots i_{V_{b}}(\om \wedge \eta )
 =
 i_{V_{1}} \cdots i_{V_{b}}(\om \wedge \eta').
 $
Since $\om$ is a basic form, we have
$
i_{V_{1}} \cdots i_{V_{b}}(\om \wedge \eta )
=
(-1)^{\ell b}\om \wedge  (i_{V_{1}} \cdots i_{V_{b}}\eta )
=
(-1)^{\ell b}\om \wedge  (i_{V_{1}} \cdots i_{V_{b}}\eta' )
=
i_{V_{1}} \cdots i_{V_{b}}(\om \wedge \eta').
$

\medskip

\nt {\em \underline{Second step}. The condition (b)}.

\smallskip

 For degree reasons it suffices to prove that $i_{V_1} \cdots i_{V_b} 
(\om \wedge d\eta)= 0$. Since $\om$ is a basic form, we can write
$i_{V_1} \cdots i_{V_b} 
(\om \wedge d\eta )=  (-1)^{\ell b} \om \wedge i_{V_1} \cdots i_{V_b}d\eta 
=0.$

\medskip

\nt {\em \underline{Third step}. The condition (c).}

\smallskip

It suffices to prove that
$
\Int{R_\mathcal{F}\times [0,1[^p} \gamma < \infty.
$
where
$\gamma \in 
\lau{\Pi}{m+p}{\mathcal{F}  \times \mathcal{I}}{M \times [0,1[^p}$ with 
compact support.
We proceed by induction on the depth of $\SF$. When the foliation is 
regular then the result is clear. In the general case we know that the 
result is true for $M-S_{_{min}}\times [0,1[^p$ and $\left(T_{_{min}} - 
S_{_{min}}\right) \times [0,1[^p$. It 
remains to consider $T_{_{min}}\times [0,1[^p$.
From  3.4.1 (d) we 
know that we can identify the perverse forms of
 $T_{_{min}} \times [0,1[^p $
 with the perverse forms of
 $D_{_{min}} \times [0,1[^{p+1}$ through the map
$$
\nabla_{_{min}} \times \Ide_{[0,1[^p}
\colon
 D_{_{min}} \times [0,1[ \times [0,1[^p
\equiv D_{_{min}} \times [0,1[^{p+1}\TO T_{_{min}}  \times [0,1[^p.
$$
 Since this map is a diffeomorphism
between
$D_{_{min}} \times ]0,1[ \times [0,1[^p$
and
$\left(T_{_{min}} -S_{_{min}}\right) \times [0,1[^p$, then we have
$$
\int_{R_{\mathcal{F}_{T_{_{min}}}} \times [0,1[^p}\gamma =
\int_{R_{\mathcal{F}_{D_{_{min}} }\times ]0,1[ \times [0,1[^{p}}} \gamma
=
\int_{R_{\mathcal{F}_{D_{_{min}} }\times [0,1[ \times [0,1[^{p}}} \gamma
=
\int_{R_{\mathcal{F}_{D_{_{min}} }\times [0,1[^{p+1}}} \gamma.
$$
The induction hypothesis gives that 
this integral is finite.

\medskip

\nt {\em \underline{Fourth step}. The condition (d).}

\smallskip

Since $\supp \om$ is compact then it suffices to 
prove 
${\displaystyle \int_{U \cap R_{\mathcal{F}}} d(\om \wedge \eta) =0}$ 
where $(U,\phii)$ is a conical chart of $\mathcal{F}$ and $\om \in 
\lau{\Om}{\ell-1}{\per{t},c}{U/\mathcal{F}}$ with $\supp \om \subset U$. 
We have:
\begin{eqnarray*}
\int_{U \cap R_{\mathcal{F}}} d(\om \wedge \eta) &=&
\int_{\R^{m-n-1} \times R_{\mathcal{G}}\times ]0,1[} d(P_{\phibas}^{*}\om 
\wedge P_{\phibas}^{*}\eta) 
\stackrel{\refp{pfi}}{=} 
\int_{\R^{m-n-1} \times R_{\mathcal{G}}\times [0,1[} 
d(\om_{\phii} 
\wedge\eta_{\phii}) \\
&\stackrel{Stokes}{=\!\!=\!\!=\!\!=} &
\int_{\R^{m-n-1} \times R_{\mathcal{G}}\times \{0 \} }\om_{\phii} 
\wedge\eta_{\phii} 
\stackrel{3.4.1 (c)}{=\!\!=\!\!=\!\!=} 0.
\end{eqnarray*}
This ends the proof.
\qed

\prgg {\bf The pairing}. Let $\eta$ be a tangent volume form.
Consider the differential operator given by 
$$
P_M \colon
\lau{\Om}{*}{\per{p}}{\mf} \times 
\lau{\Om}{\ell-*}{\per{q},c}{\mf} \TO \R
$$
where
$
P_M (\alpha,\beta) = {\displaystyle \int_{R_\mathcal{F}} }\alpha \wedge \beta \wedge \eta.
$
Notice that the manifold $R_{\F} \subset M$ is an 
oriented manifold. 
From Proposition \ref{etatang} (c) we have that this integration is 
well-defined.
This operator depends on the action $\Phi \colon G \times M\to M$ and 
on the choice of the basis $\{u_{1}, \ldots ,u_{b}\}$  of 
$\mathfrak{g}$ (cf.  Proposition 
\ref{etatang} (a)).

The {\em pairing} is the induced operator
$$
P_M 
 \colon
\lau{\IH}{*}{\per{p}}{\mf} \times 
\lau{\IH}{\ell-*}{\per{q},c}{\mf} \TO \R
$$
defined by
$$
P_M ([\alpha],[\beta]) = \int_{R_\mathcal{F}} \alpha \wedge \beta \wedge \eta.
$$
This operator is well-defined (cf. Proposition \ref{etatang} (b) and 
(d)).
The PoincarŽ Duality stands that $P_M$ is a non degenerate pairing, that 
is, the operator
$$
P_M 
 \colon
\lau{\IH}{*}{\per{p}}{\mf} \TO \Hom 
\left(\lau{\IH}{\ell-*}{\per{q},c}{\mf} , \R
\right)
$$
defined by
$$
P_M ([\alpha])([\beta]) = \int_{R_\mathcal{F}} \alpha \wedge \beta \wedge \eta
$$
is an isomorphism. 

The first step to get the PoincarŽ Duality is the following.

\bP
\label{modelDP}
Consider a twisted 
product $K \times_{H }\R^n$ as in 2.2.  Suppose that the action
$\Phi \colon G \times (K \times_{H }\R^n) \to (K \times_{H }\R^n)$ is 
effective. We fix a basis of 
$\mathfrak{k}$ as in 3.8.3.
The BIC of $(K \times_{H }\R^n,\F_{tw})$
verifies the PoincarŽ Duality when the BIC of its slice  
$(\R^n,\F_{\R^n})$ verifies the PoincarŽ Duality.
\eP
\pro
Since  the action $\Phi \colon G \times (K\times_{H}\R^n) \to 
(K\times_{H}\R^n)$ defining $\F_{tw}$ is  effective then the action 
$\Theta \colon (G \cap H) \times \R^n \to \R^n$  defining
$\F_{\R^n}$ is also effective. The Proposition makes sense.
We proceed  now in two steps.

\Zati \underline{\em Construction of the tangent volume form}.

Consider a tangent volume form 
$\eta _0\in \left(\lau{\Pi}{b-a}{\mathcal{F}_{\R^n}}{\R^n}\right)^H$ 
  of the slice $(\R^n,\mathcal{F}_{\R^n})$. 
 We construct the tangent volume form $\eta$ of the twisted 
product in terms of $\eta_{0}$.

Put $\eta = \nabla^{-1} 
(\gamma_1 \wedge \cdots \wedge \gamma_a \wedge \eta_0)=
R^{-*}(\gamma_1 \wedge \cdots \wedge \gamma_a \wedge \Delta'(\eta_0)).
$
It belongs to $
\left(\lau{R}{*}{\mathcal{F}}{K \times_H \R^n}\right)^K 
$
(cf. 3.8.5).
We prove 
now \refp{eta}.
Since $R_*X_i = V_i$ for $i \in \{ 1, \ldots, f\}$ then,
$$
\begin{array}{rcl}
\eta (V_1,\ldots,V_b) &= &
(\gamma_1\wedge \cdots \wedge \gamma_a \wedge
\Delta'(\eta_0))(X_1,\ldots,X_b) = \Delta'(\eta_0)(X_{a+1},\ldots,X_b) 
=\\[,2cm]
&=&
 \left(\eta_0 + \suma{a+1 \leq i_1 < \cdots < i_\ell \leq c}{} \ \gamma_{i_1} 
\wedge \cdots 
\wedge \gamma_{i_\ell} \wedge 
\left(i_{W{_{i_\ell}} }\cdots 
i_{W{_{i_1}}}\eta_0\right)\right)(X_{a+1},\ldots,X_b)= \\[,6cm]
&=&i_{W{_{b}} }\cdots i_{W{_{a+1}}}\eta_0 = 1
\end{array}
$$
from $\refp{eta}$ for $\eta_{0}$. We conclude that $\eta$ is a tangent 
volume form of the twisted product.
Recall that the pairing $P_{M}$ can be defined using this form $\eta$ 
(cf. Proposition \ref{etatang} (a)).

\Zati \underline{\em PoincarŽ Duality}.

Consider now two complementary perversities $\per{p}$ and $\per{q}$ on $\dim 
K \times_H \R^n$. The induced perversities $\per{p}$ and $\per{q}$ on 
$\R^n$ are 
also two complementary perversities. Let us see that, for each stratum
$S\in \SFRn$, we have:
\begin{eqnarray*}
\per{p} (S)+ \per{q}(S)
&=&
\per{p} (K \times_H S)+ \per{q}(K \times_H S)
=
\per{t}(K \times_H S)  =\\
&=&
\dim K \times_H \R^n - \dim K \times_H S - \dim G_{K \times_H S}-2
 =\\
&=& \dim \R^n - \dim S - \dim (G\cap H)_S -2
=
\per{t}(S) .
\end{eqnarray*}
By 
hypothesis  the pairing
$
P_{\R^n} \colon \lau{\IH}{*}{\per{p}}{\R^n/\mathcal{F}_{\R^n}}\times 
\lau{\IH}{\ell_{\R^n}-*}{\per{q},c}{ \R^n/\mathcal{F}_{\R^n}} \to \R
$
is non degenerate. 
Since  $\eta_{0}$ is 
$H$-invariant then the pairing
$
P_{\R^n} \colon \left( 
\lau{\IH}{*}{\per{p}}{\R^n/\mathcal{F}_{\R^n}}\right)^{H}\times 
\left(\lau{\IH}{\ell_{\R^n}-*}{\per{q},c}{ 
\R^n/\mathcal{F}_{\R^n}}\right)^{H} \to \R
$
is also non degenerate. On the other hand, it is clear that the pairing
$
P \colon 
{\bigwedge}^*(\gamma_{c+1},\ldots,\gamma_{f}) 
\times 
{\bigwedge}^{f-c-*}(\gamma_{c+1},\ldots,\gamma_{f})
\TO \R,
$
defined by 
$$
P(\gamma_{i_1} \wedge \cdots \wedge \gamma_{i_u},
\gamma_{j_1} \wedge \cdots \wedge \gamma_{j_{f-c-u}})
=
{\displaystyle \int_{K/H}  \gamma_{i_1} \wedge \cdots \wedge \gamma_{i_u}}
\wedge \gamma_{j_1} \wedge \cdots \wedge \gamma_{j_{f-c-u}}
$$
is non degenerate. Notice the equality $\ell_{K\times_H \R^n}    = f+a + n -c-b = 
\ell_{\R^n} + f-c$.

We prove that the pairing 
$
P_{K \times_H \R^n} \colon \lau{\IH}{*}{\per{p}}{K \times_H \R^n}\times 
\lau{\IH}{\ell_{K\times_H \R^n}  - *}{\per{q},c}{K \times_H \R^n} \TO \R
$
is not degenerate.
We know, from Proposition \ref{isom}, that this is the case if the following 
diagram commutes (up to a sign):
$$
\begin{CD}
\left({\bigwedge}^*(\gamma_{c+1},\ldots,\gamma_{f})
\otimes \left( 
\lau{\IH}{*}{\per{p}}{\R^n/\mathcal{F}_{\R^n}}\right)^H\right)
\! \times \!
\left({\bigwedge}^{f-c-*}(\gamma_{c+1},\ldots,\gamma_{f})
\otimes \left( 
\lau{\IH}{\ell_{\R^n}-*}{\per{p}}{\R^n/\mathcal{F}_{\R^n}}\right)^H\right)
@>P \otimes P_{\R^n}>> \R\\
    @V \nabla^{-*} \times \nabla^{-*}  VV @V\Ide VV\\
\lau{\IH}{*}{\per{p}}{K \times_H \R^n}\times 
\lau{\IH}{\ell_{K\times_H \R^n} -*}{\per{q},c}{K \times_H \R^n}
@>P_{K \times_H \R^n} >> \R\\
    \end{CD}
    $$
Let us see that. For each $[\alpha] \in 
\left(\lau{\IH}{i}{\per{p}}{\R^n/\mathcal{F}_{\R^n}}\right)^H$
and 
each $[\beta] \in \left(\lau{\IH}{\ell_0-i}{\per{p}}{\R^n/\mathcal{F}_{\R^n}}\right)^H$, 
we have, for degree reasons,
$$
\begin{array}{l}
P_{K\times_H \R^n} (\nabla^{-1}\times \nabla^{-1}) (\gamma_{i_1} \wedge \cdots \wedge \gamma_{i_u}
\otimes [\alpha],\gamma_{j_1} \wedge \cdots \wedge \gamma_{j_{f-c-u}}
\otimes [\beta])
= 
\\[,3cm]
\Int{K\times_H \R^n}R^{-*}
(\gamma_{i_1} \wedge \cdots \wedge \gamma_{i_u} \wedge \alpha  
\wedge
\gamma_{j_1} \wedge \cdots \wedge \gamma_{j_{f-c-u}}\wedge \beta 
\wedge
\gamma_1 \wedge \cdots \wedge \gamma_a \wedge \eta_0)=
\\[,3cm]
\Int{K\times \R^n}
\gamma_{i_1} \wedge \cdots \wedge \gamma_{i_u} \wedge \alpha  
\wedge
\gamma_{j_1} \wedge \cdots \wedge \gamma_{j_{f-c-u}}\wedge \beta 
\wedge
\gamma_1 \wedge \cdots \wedge \gamma_a \wedge \eta_0 \wedge \gamma_{a+1} 
\wedge \cdots \wedge \gamma_c=
\\[,3cm]
\Int{K}
\gamma_1 \wedge \cdots \wedge \gamma_c \wedge
\gamma_{i_1} \wedge \cdots \wedge \gamma_{i_u} 
\wedge
\gamma_{j_1} \wedge \cdots \wedge \gamma_{j_{f-c-u}}
\cdot 
\Int{\R^n}
\alpha  
\wedge
\beta 
\wedge
\eta_0=
\\[,3cm]
 P(\gamma_{i_1} \wedge \cdots \wedge \gamma_{i_u} ,
\gamma_{j_1} \wedge \cdots 
\wedge \gamma_{j_{f-c-u}}) \cdot
P_{\R^n}([\alpha],[\beta])=
\\[,3cm]
(P \otimes P_{\R^n})(\gamma_{i_1} \wedge \cdots \wedge 
\gamma_{i_u}\otimes [\alpha],\gamma_{j_1} \wedge \cdots 
\wedge \gamma_{j_{f-c-u}} \otimes [\beta]),
\end{array}
$$
up to a sign.
\qed

We study now the PoincarŽ Duality of the tubular neighborhoods of the 
strata of $\SF$ (cf. 2.3).
Consider $(T,\tau,S,\R^n)$  a $K$-invariant tubular neighborhood of a 
singular stratum $S$, where $K$ is a tamer group of $G$ with $K = 
\overline{G}$.
Put  $(\R^n,\F_{\R^n})$ the slice of the tubular neighborhood. The foliation 
$\F_{\R^n}$ is defined by an effective tame action $\Theta \colon 
(G_{S})_{0} \times \R^n \to \R^n$.
We fix $\{u_{1}, \ldots , u_{b}\}$ 
a basis of $\mathfrak{g}$ and suppose that $\{u_{a+1}, \ldots u_{b}\}$ 
is a basis of of the Lie algebra of $G_{S}$.

\bP
\label{tubDP}
Under the above conditions, 
if the BIC of the slice $(\R^n,\mathcal{F}_{\R^n})$ verifies the PoincarŽ Duality then 
the BIC of the tube
$(T,\mathcal{F})$ also verifies the PoincarŽ Duality
\eP
\pro The proof of the Proposition is the same of that of the Proposition \ref{tub} 
by
changing

\bigskip

	\begin{tabular}{lcl}
$\mathfrak{A}(T,\mathcal{F})$& $=$&
``The BIC
	$
	\lau{\IH}{*}{\per{p}}{T/\mathcal{F}}
	$
	is finite dimensional for each perversity $\per{p}$."\\[,3cm]
 by &&\\[,3cm]
$\mathfrak{A}(T,\mathcal{F}) $& $=$&
``The pairing
	$
	P_{T} \colon \lau{\IH}{*}{\per{p}}{T/\mathcal{F}} \TO
\lau{\IH}{\ell - *}{\per{q}}{T/\mathcal{F}} 
	$
	is non degenerate, \\
	&& for any two complementary perversities $\per{p}$ and $\per{q}$."
\end{tabular}

\medskip

We consider the orbit type stratification 
of $S$  induced by the action $\Phi_{S} 
\colon K \times S \to S$ of $\Phi$. We proceed by induction on the depth of 
this stratification.
By using the Mayer-Vietoris technics of 3.6 and 3.7
we can 
suppose that $\Phi_{S}$ defines  a fiber bundle $\pi \colon S \to 
S/K$ whose fiber is $K/H$. Considering a good covering of $S/K$, the 
Mayer-Vietoris procedure leads us to the case $S/K = point$, that is, 
to the case where $T$ is the twisted product $K\times_H \R^n$. Here, we apply the Proposition 
\ref{modelDP}.\qed

The second main result of this section is the following 
\bT
\label{TPC}
The BIC of the foliation determined by an isometric action of an abelian 
Lie group on an oriented 
compact manifold verifies the PoincarŽ Duality.
\eT
\pro In fact, we prove the result for a foliation $\F$ determined by 
a tame action $\Phi \colon G \times M \to M$ on an oriented  manifold not 
necessarily compact. We can suppose that $G$ is connected and that the 
action $\Phi$ is an effective one. We fix $K$ a tamer group where $G$ 
is dense.

The proof of the Proposition is the same of that of the Theorem \ref{T} by
changing

\bigskip

	\begin{tabular}{lcl}
$\mathfrak{A}(M,\mathcal{F})$& $=$&
``The BIC
	$
	\lau{\IH}{*}{\per{p}}{M/\mathcal{F}}
	$
	is finite dimensional for each perversity $\per{p}$."\\[,3cm]
 by &&\\[,3cm]
$\mathfrak{A}(M,\mathcal{F}) $& $=$&
``The pairing
	$
	P_{T} \colon \lau{\IH}{*}{\per{p}}{M/\mathcal{F}} \TO
\lau{\IH}{\ell - *}{\per{q}}{M/\mathcal{F}} 
	$
	is non degenerate, \\
	&& for any two complementary perversities $\per{p}$ and $\per{q}$."
\end{tabular}

\smallskip

 By using the Mayer-Vietors technics of 3.6, 3.7  and the Proposition 
\ref{tubDP} we reduce the problem to prove 
$
\mathfrak{A}(\R^n,\mathcal{F}_{\R^n}) ,
$
where 
$(\R^n,\mathcal{F}_{\R^n}) 
\equiv
\left( c\S^{n -1} , 
c\mathcal{G}_{S}\right)$
is a slice of a tubular neighborhood of a singular stratum $S$ of $\SF$. In other words, we need to prove that the 
pairing
$$
P_{\R^n} \colon \lau{\IH}{*}{\per{p}}{\R^n/\mathcal{F}_{\R^n}}  \times
\lau{\IH}{\ell_{\R^n} -*}{\per{q},c}{\R^n/\mathcal{F}_{\R^n}}  \TO \R,
$$
is non degenerate.

From Propositions \ref{pro}, \ref{local5}, \ref{procomp} and \ref{local6} 
we have
\begin{eqnarray}
\label{cono1}
\lau{\IH}{i}{\per{p}}{\R^n/\mathcal{F}_{\R^n}} &=&
\left\{
\begin{array}{cl}
\lau{\IH}{i}{\per{p}}{\S^{n-1}/\mathcal{G}_{S}}  & \hbox{if } i \leq 
\per{p}(\vartheta)\\[,3cm]
0 & \hbox{if } i \geq 
\per{p}(\vartheta)+1
\end{array}
\right.
\end{eqnarray}
and
\begin{eqnarray*}
\lau{\IH}{\ell_{\R^n} -i}{\per{q},c}{\R^n/\mathcal{F}_{\R^n}} &=&
\left\{
\begin{array}{cl}
0 & \hbox{if } i\geq 
\ell_{\R^n} - \per{q}(\{ \vartheta\}) - 1 \\[,3cm]
\lau{\IH}{\ell_{\R^n} -i-1}{\per{q}}{\S^{n-1}/\mathcal{G}_{S}}  & 
\hbox{if } i \leq 
\ell_{\R^n}  - \per{q}(\{ \vartheta\}) - 2.
\end{array}
\right.
\end{eqnarray*}
Since $\per{p}$ and $\per{q}$ are complementary perversities on $\R^n$ 
then we have:
$$
\per{p}(\{ \vartheta\}) + \per{q}(\{ \vartheta\})
=
\per{t}(\{ \vartheta\}) 
=
n -\dim\mathcal{F}_{\R^{n}} - 2 =n - 
\dim\mathcal{G}_{S} -2 = \ell_{\S^{n-1}} -1 = \ell_{\R^n} -2.
$$
These 
formul\ae\ give:
\begin{equation}
\label{cono2}
\lau{\IH}{ \ell_{\R^{n}}  - i}{\per{q},c}{\R^n/\mathcal{F}_{\R^n}} =
\left\{
\begin{array}{cl}
0 & \hbox{if } i \geq 
\per{p}(\{ \vartheta\}) + 1\\[,3cm]
\lau{\IH}{\ell_{\S^{n-1}}-i}{\per{q}}{\S^{n-1}/\mathcal{G}_{S}}  & \hbox{if } 
\leq 
\per{p}(\{ \vartheta\}).
\end{array}
\right.
\end{equation}
Now, an argument on the depth of $\SF$ gives 
the property $\mathfrak{A}(\R^n,\mathcal{F}_{\R^n})$ from these three facts
\begin{itemize}
\item[(i)] $\mathfrak{A}(\S^{n-1},\mathcal{G}_{S})$.

\item[(ii)] The pairing $P_{\R^n}$ becomes the 
pairing $P_{\S^{n-1}}$ through the
isomorphism induced by \refp{cono1} and \refp{cono2}.
\begin{quote}
    The operator
     $
     \aleph_{23} \colon 
     \lau{\IH}{*}{\per{p}}{\S^{n-1}/\mathcal{G}_S}
     \to
     \lau{\IH}{*}{\per{p}}{c\S^{n-1}/
     c\mathcal{G}_S}
     $  
     defining \refp{cono1} is 
     $
     \aleph_{23}([\alpha]) = [\alpha];
     $  the operator 
     $
     \aleph_{24} \colon 
     \lau{\IH}{*}{\per{q}}{\S^{n-1}/\mathcal{G}_S}
     \to
     \lau{\IH}{*}{\per{q},c}{c\S^{n-1}/
     \mathcal{G}_S}
     $  
     defining \refp{cono2} is $
     \aleph_{24}([\beta]) = [g\, 
     dt \wedge \beta].
     $ 
Since the action of $\Theta$ lies on $\S^{n-1}$ then 
we can take a common tangent volume form $\eta$ for 
$\mathcal{F}_{\R^n}$ and $\mathcal{G}_S$ (cf. proof of Proposition 
\ref{etatang}). Now, for 
$[\alpha] \in \lau{\IH}{i}{\per{p}}{\S^{n-1}/\mathcal{G}_S}$
and
$
[\beta] \in
\lau{\IH}{\ell_{\S^{n-1}}-i}{\per{q},c}{\S^{n-1}/\mathcal{G}_S}
$
we have
\begin{eqnarray*}
P_{\R^n} (\aleph_{23}[\alpha], \aleph_{24}[\beta])&=&
\int_{\R^{n} \times R_{\mathcal{F}_{S^{n-1}}} \times ]0,1[} \alpha 
\wedge  g \, 
\wedge dt \wedge \beta \wedge \eta =\\
&=&
\left(\int_{ R_{\mathcal{F}_{S^{n-1}}}} \alpha \wedge \beta \wedge \eta\right)
\left(\int_{0}^1 g dt\right)
= 
P_{\S^{n-1}} ([\alpha], [\beta]).
\end{eqnarray*}

\end{quote}
\item[(iii)] The perversities $\per{p}$ and $\per{q}$ are complementary on 
$\S^{n-1}$.
\begin{quote}
We have, for any stratum $S \in \hbox{\sf S}_{{\mathcal G}}$ 
the equalities
$\per{p}(S) + \per{q}(S) =\per{p}(S \times ]0,1[) + \per{q}( S \times ]0,1[)=
\per{t}(S \times ]0,1[) 
= \codim_{\R^n} \mathcal{F}_{\R^n} - 
\codim_{S \times ]0,1[} (
\mathcal{F}_S \times \mathcal{I}) -2 = \codim_{\S^{n-1}} 
\mathcal{G}_S -\codim_S \mathcal{F}_S -2 = \per{t}(S)$.
\end{quote}
\end{itemize}
Hau amaia da. \qed


%

\end{document}